\documentclass[11pt,
]{article}

\usepackage[tbtags]{amsmath}
\usepackage{amsfonts,amssymb,mathrsfs,amscd,amstext,
latexsym,amsbsy,amsthm}

\usepackage{bbold}
\usepackage{mparhack}
\usepackage{marginnote}

\usepackage{imakeidx}
\indexsetup{level=\subsection*,toclevel=subsection}

\input{36.sty}

\makeindex

\begin{document}

\title
{Definable minimal collapse functions at arbitrary
projective levels
}

\author 
{Vladimir Kanovei\thanks{IITP RAS and MIIT, Moscow, Russia,
 \ {\tt kanovei@googlemail.com} --- 
contact author.
Partial support of RFFI grant 03-01-00757 acknowledged.}  
\and
Vassily Lyubetsky\thanks{IITP RAS,  Moscow, Russia, 
\ {\tt lyubetsk@iitp.ru}.
Supported by RSF grant 14-50-00150.
}}

\date 
{\today}

\maketitle


\begin{abstract}
Using a non-Laver modification of Uri Abraham's minimal
$\varDelta^1_3$ collapse function,
we define a generic extension $L[a]$ by a real $a$,
in which, for a given $n\ge3$,
$\{a\}$ is a lightface $\varPi^1_n$ singleton,
$a$ effectively codes a cofinal map $\omega\to\omega_1^L$
minimal over $L$,
while every $\varSigma^1_n$
set $X\subseteq\omega$ is still constructible.
\end{abstract}

\markright{Definable collapse functions
at arbitrary projective levels}

{\small
\def\contentsname{\large Contents}
\tableofcontents
}

\parf{Introduction}
\las{or}

It is well-known that all sets $x\sq\om$ of   
the lightface class $\is12$ or $\ip12$ 
are Goedel-constructible. 
In fact this is an immediate corollary of the 
Shoenfield absoluteness theorem. 
But one gets models with nonconstructible sets 
which belong to the analytic hierarchy just above 
the mentioned threshold. 
In particular it is consistent with $\ZFC$ 
that there exists a $\id13$ real and a $\ip12$ 
real singleton, see \cite{jsad}, and such a real 
can be of minimal \dd\rL degree, \cite{jenmin}. 

Many more results on definable sets of different 
kind have been obtained on the base of forcing 
methods invented in the abovementioned papers. 
Most of them employ versions of the 
\rit{almost disjoint} coding method of \cite{jsad}. 
A recent article \cite{ffmm} contains several 
powerful applications of almost disjoint coding, 
in particular, to the construction of 
models with $\id13$ well-orderings of the reals, in  
which the reals have some very special properties.
The paper also contains further references. 

Yet the almost disjoint coding technique is pretty 
useless in the case of models containing definable 
generic objects and \rit{minimal} over the 
ground model with respect to this or another property. 
The first example of such a model was presented by 
Jensen \cite{jenmin}. 
Namely, Jensen's forcing notion $\jf\in\rL$ consists 
of perfect trees in $\nse$ 
(a subset of the Sacks forcing), 
and if a real $a\in\bn$ is \dd\jf generic over $\rL$ 
then 1) it is true in $\rL[a]$ that $\ans a$ is a 
nonconstructible $\ip12$ singleton, and 2) $a$ is 
\rit{minimal over} $\rL$, in the sense that if 
$b\in\rL[a]\cap\bn$ then either $b\in\rL$ or 
$a\in\rL[b]$. 
(See also 28A in \cite{jechmill} on this forcing.)

Several variations of this forcing are known. 
In particular, a model in \cite{ian78} in which, 
for a given $n\ge3$ there exists a minimal 
nonconstructible $\ip1n$ singleton but all $\is1n$ 
sets $x\sq\om$ are constructible, 
an \dd{\om_2}long 
iteration of Jensen's forcing in \cite{abr}, 
or a recent model 
in \cite{kl22} in which there is 
an equivalence class of the equivalence
\index{equivalence relation!E0@$\Eo$}%
\index{zzE0@$\Eo$}%
relation $\Eo$\snos
{\label{deo}%
$\Eo$ is defined on the Baire space $\bn$ 
so that $x\Eo y$ iff the set $\ens{n}{x(n)\ne y(n)}$
is finite.} 
(a \dd\Eo\rit{class}, for brevity), 
which is a lightface $\ip12$ set in $\bn$,
not containing  $\od$ elements, and a related 
model in \cite{kl25} containing a $\ip12$ 
Groszek -- Laver pair of \dd\Eo classes. 

The research of this paper was inspired by another 
minimal-style forcing construction, a generic 
extension $\rL[a]$ by Abraham~\cite{abr2} such that  
1) $\ans a$ is a 
nonconstructible $\ip12$ singleton in $\rL[a]$, 
2) $\omi^{\rL[a]}=\om_2^{\rL}$ (so $a$ codes a 
collapse of $\omi^\rL$), 
and 
3) $a$ is \rit{a minimal collapse over} $\rL$, 
in the sense that if 
$b\in\rL[a]$, $b:\om\to\omi^\rL$, and 
$b$ is cofinal in $\omi^{\rL}$, 
then $a\in\rL[b]$. 
Abraham's forcing in \cite{abr2} consists of 
Laver-style trees in $\bsd$, and its complicated 
construction in $\rL$, 
while having a certain semblance of 
Jensen's method in \cite{jenmin}, 
involves some crucial novel ideas.

Our main result extends this research line.  
The next theorem asserts the existence of
a model of $\zfc$, in which, for a given $\nn\ge2$, 
there is a $\ip1\nn$ real singleton which codes 
the collapse of $\omi^{\rL}$ in minimal way, 
and in the same time reals in $\is1\nn$ do not 
code the collapse.
The abovementioned result of \cite{abr2} corresponds to
the case $\nn=2$ in this theorem, of course. 
We use the blackboard $\nn$ to distinguish the fixed 
number $\nn$ in the theorem from other numbers $n$
in the text.

\bte
\lam{Tun}
Let\/ $\nn\ge2$. 
There is a generic extension\/ $\rL[a]$ of\/ $\rL$, 
by a real\/ $a\in\bn$, such that the following is 
true in\/ $\rL[a]$$:$ \ \ \ \ \ 
{\rm(i)}  $\omi^{\rL[a]}=\om_2^{\rL}\,;$
%
\ben
\renu
\atc
\itlb{Tun2}%
{\rm(minimality)}
if\/ $b\in\rL[a]$, $b:\om\to\omi^\rL$,   
$b$ is cofinal in $\omi^{\rL}$, 
then\/ $a\in\rL[b]\,;$ 

\itlb{Tun3}%
the singleton\/ $\ans{a}$ is a 
(lightface)\/ $\ip1\nn$ set$;$

\itlb{Tun4}%
{\rm(vacuous for $\nn=2$)}
every\/ $\is1{\nn}$ set\/ $x\sq\om$ 
belongs to\/ $\rL$.
\een
\ete

\parf{Structure of the proof}
\las{ko8}

The proof of Theorem~\ref{Tun} is organized as follows. 
Basic notions, related to \dd\omi branching trees 
in $\bsd$ (\rit{wide trees\/}), 
are introduced in Section~\ref{tre}. 
Unlike \cite{abr2}, we'll not focus on 
Laver-style trees, which makes basic constructions 
somewhat simpler.
Every set ${\jpi}$ of wide trees $T$, closed under 
restrictions, 
is considered as \rit{a forcing by wide trees}, 
a \rit{\sfo} in brief, Section~\ref{wds}. 
Every \sfo\ adjoins a \dd{\jpi}generic ``real'' 
$a\in\on.$ 

Section~\ref{ite} presents a non-Laver  
modification of Abraham's method in 
\cite[2.14]{abr2}, 
designed to define uncountable decreasing sequences 
of wide trees.
Basically, any collection $F$ of wide trees, satisfying 
rather transparent conditions of Definition~\ref{dri}, 
yields a wide tree $\tre F$, so that if $F\sq F'$ 
then $\tre{F'}\sq \tre F$. 
We apply the method to prove Theorem~\ref{frt} in 
sections \ref{nef}, \ref{pret}, which allows, given a 
wide tree $S$ and a family of continuous functions 
$f_\al:\on\to\on$, $\al<\omil$, 
to define a smaller wide tree $T\sq S$, 
regular in some sense with respect to each $f_\al$. 

Another technical device, also having its roots in 
\cite{abr2}, is introduced in Section~\ref{bc}. 
It allows to shrink a given wide tree $S$ to a smaller wide 
tree $T$ such that any pre-dense set $U\sq S$ in a 
given family of \dd\ali many such sets
meets every infinite branch in $T$ except for a bounded 
set of them 
(Corollary~\ref{plo3}).

Then, arguing in the constructible universe $\rL$, we  
define a forcing notion to prove Theorem~\ref{Tun} 
in Section~\ref{vop} 
in the form $\dP=\bigcup_{\al<\omd}\dP_\al$. 
The summands $\dP_\al$ are \dd\ali large \sfo s 
defined by induction. 
Any \dd\dP generic extension of $\rL$ happens to be a 
model for Theorem~\ref{Tun}, which we prove in 
the remainder. 

The inductive construction of $\dP_\al$ involves two 
key genericity ideas.
The first idea, essentially by Jensen~\cite{jenmin}, 
is to make every level $\dP_\al$ of the construction  
\rit{generic} in some sense over the union of lower  
levels $\dP_\xi$, $\xi<\al$. 
This is based on a construction developed in  
sections~\ref{wtfn}, \ref{vop}, which includes   
the abovementioned modification of Abraham's method.
The iterated genericity of the 
levels $\dP_\al$ implies that the two sets are 
equal in any \dd\dP generic extension of $\rL$:\vom 

1) the singleton $\ans{\xx G}$ of the principal
generic element $\xx G\in\on$,\vom 

2) the intersection 
$\bigcap_{\al<\omd}\bigcup_{T\in\dP_\al}$.\vom 

\noi
This equality, eventually leading to \ref{Tun2} of 
Theorem~\ref{Tun}, is established in 
sections \ref{pres+},
\ref{bex1}, 
on the base of studies of continuous functions 
in sections \ref{nef}, \ref{pret}. 

The second idea goes back to old papers 
\cite{h74}, \cite{ian78}.
In $\rL$, let $\vstf$ be the set of all countable 
sequences $\vjpi=\sis{\jpi_\xi}{\xi<\al}$ 
($\al<\omi$), 
compatible with the first genericity idea at each
step $\xi<\al$. 
Then a whole sequence $\sis{\jpi_\al}{\al<\omi}$ 
can be interpreted as a maximal chain in $\vstf$. 
It happens that if such a chain is \rit{generic}, in 
some sense precisely defined in Section~\ref{vop}, 
\ref{ep2} of Theorem~\ref{ep},  
with respect to all $\fs1{\nn-1}$ subsets of $\vstf$, 
then the ensuing forcing notion 
$\dP=\bigcup_{\al<\omi}\dP_\al$ 
inherits some basic forcing properties of the whole 
forcing by (all) wide trees, up to the \dd\nn th level 
of projective hierarchy. 
This includes, in particular, the invariance of the 
forcing relation with respect to some natural 
transformations of wide trees, leading eventually 
to the proof of \ref{Tun4} of Theorem~\ref{Tun} 
in sections \ref{fs} -- \ref{finar}.

\parf{Wide trees} 
\las{tre}

Let $\bsd$ be the set of all \rit{strings}   
\index{zzom1om@$\bsd$}%
\index{string!$\bsd,$ string}%
(finite sequences) 
\index{string}%
of ordinals $\xi<\omi$ --- 
including  \rit{the empty string} $\La$.
\index{zzLa@$\La$, the empty string}%
\index{string!$\La$, the empty string}%
\index{string!extension, $\sq$, $\su$}%
If $s\in\bsd$ then $\lh s<\om$ is the \rit{length}
of a string $s$, and 
$\tmax s<\omi$ is the largest term in $s$.
\index{zzlhs@$\lh s$, length}%
\index{length!$\lh s$}%
\index{zzmaxs@$\tmax s$, largest term}%
\index{string!length}%
\index{string!$\omi^n,$ strings of length $n$}%
Let 
$\omi^n=\ens{s\in\bsd}{\lh s=n}$
(strings of length $n$). 
If $t\in\bsd$ and $\xi<\omi$, then $t\we\xi$ is
the extension of $t$ by $\xi$ as the rightmost 
\index{zztiw@$t\we\xi$}%
term. 
If $s,t\in\bsd$ then $s\sq t$ means that the string $t$ 
\index{zzstsq@$s\sq t$}%
\index{zzstsu@$s\su t$}%
\index{string!extension, $\sq$, $\su$}%
extends $s$, while 
$s\su t$ means a proper extension. 
A set $T\sq\bsd$ is a \rit{tree}
\index{tree}
iff $s\in T\imp t\in T$ 
whenever $s,t\in\bsd$ and $t\su s$. 
Then: 
\bit
\item[$-$]
if $s\in T$ then   
$\suc Ts=\ens{t\in T}{s\su t\land \lh t=\lh s+1}$,
\imar{suc Ts}%
\index{successor!$\suc Ts$}%
\index{zzsuccTs@$\suc Ts$}%
the set of all \rit{successors} of $s$ in $T$. 
If $\suc Ts=\pu$ then $s$ is an \rit{\eno} 
of $T$; 

\item[$-$]
$\br T=\ens{s\in T}{\card{(\suc Ts)}\ge 2}$,
\imar{br T}%
all \rit{\bno s} of $T$, and \\
$\brn Tn=\ens{s\in\br T}
{\card{(\ens{u\in\br T}{u\su s})}=n}$;

\item[$-$]
if $u\in T$ then define 
\index{tree!restricted $T\ret u$}%
\index{zzTIu@$T\ret u$}%
$T\ret u=\ens{t\in T}{u\sq t\lor t\sq u}$,
a \rit{restricted tree}; 

\item[$-$]
if $T$ is not pairwise \dd\sq compatible then there is 
\index{tree!stem, $\roo T$}%
\index{stem, $\roo T$}%
\index{zzstemT@$\roo T$}%
a largest string $u\in T$ such that $T\ret u=T$,
denoted by $u=\roo T$, then 
$\ans{\roo T}=\brn T0$;

\item[$-$]
$[T]=\ens{x\in\on}{\kaz m\,(x\res m\in T)}$, 
\index{zzTII@$[T]$}%
\index{tree!branch}%
a closed set in $\on$.
\eit
\vyk{
If $s\in\bsd$ then put $\bad s=\ens{x\in\on}{s\su x}$, 
\index{zzBs@$\bad s$}%
a basic clopen subset of $\on$.
}

\bdf
\lam{dent}
A set\/ $U\sq T$ is \rit{dense} in a tree $T$  
\index{set!dense}%
\index{dense}%
if $\,\kaz s\in T\,\sus u\in U\,(s\sq u)$,  
\index{set!open dense}%
\index{open dense}%
\rit{open dense}, if in addition   
$s\in U$ holds whenever $s\in T$, $u\in U$, 
$u\sq s$,  
and 
\index{set!pre-dense}%
\index{pre-dense}%
\rit{pre-dense}, if the set 
$U'=\ens{s\in T}{\sus u\in U(u\sq s)}$ is dense.
\edf

\bdf
\lam{std}
A tree $\pu\ne T\sq\bsd$ is a \rit{wide tree},
\index{tree!wide, $\pew$}%
\index{zzWT@$\pew$}%
in symbol $T\in\pew$,
if 
any $s\in T$ can be extended to a \bno\ 
$t\in\br T$, $s\sq t$, 
and if $t\in \br T$ then
$\card{(\suc Ts)}=\ali$ 
--- \ie, all \bno s 
are \dd\omi branching. 
\edf

A bigger set $\pew'$ consists of all trees
\index{tree!WT1@$\pew'$}%
\index{zzWT1@$\pew'$}%
$\pu\ne T\sq\bsd$ such that each subtree
of the form $\req Ts$, $s\in T$, is uncountable.
Clearly $\pew\sneq\pew'$, but $\pew$ is still
dense in $\pew'$, so that every
tree $T\in\pew'$ contains a subtree $S\in\pew$,
$S\sq T$.

Generally, $\pew$ and $\pew'$ belong to the category
of \rit{uncountably-splitting}
versions of the perfect set forcing.
Similar forcing notions of this kind, as well as
their Laver-style versions (which require every node
above the stem to be a wide-splitting node), have been
thoroughfully studied in set theoretic papers,   
see \eg\ Namba~\cite{nam71}, Bukovsky~\cite{buk76},
Abraham~\cite{abr2}, Jech~\cite[Chap.~28]{jechmill},
to mention a few.

\ble
\lam{iso}
Suppose that\/ $T\in\pew$. 
If\/ $s\in T$ then\/ $\req Ts\in\pew$. 
If\/ $x\in X\sq[T]$, $X$ is open in\/ $[T]$,  
then there is\/ $s\in T$ such that\/
$s\su x$ and\/ $\req Ts\sq X$. 
\qed
\ele

\bdf
\lam{sqn}
We introduce two notions of inclusion between trees
which partially honor the branching structure. 
If $S,T\sq\bsd$ are trees then define:
\bit
\item[$-$]
$S\sq_n T$ \ iff \ $\brn Tn\sq S\sq T$;
\index{zzSsubsetnT@$S\sq_n T$}%
\index{inclusion!SsubsetnT@$S\sq_n T$}%

\item[$-$]
$S\sq'_n T$ \ iff \ $S\sq T$ and $\brn S{n-1}=\brn T{n-1}$.\qed
\index{zzSsubset'nT@$S\sq'_n T$}%
\index{inclusion!Ssubset'nT@$S\sq'_n T$}%
\eit
\eDf

\ble
\label{2sq}
\ben
\renu
\itlb{2sq1}
The relations\/ $\sq_0$ and\/ $\sq_0'$ coincide with
just\/ $\sq\;;$
\itlb{2sq2}%
${S\sq'_{n+1} T}\imp{S\sq_n T}\imp{S\sq'_n T}\;;$ 

\itlb{2sq3}%
if $S\sq_n T$ then\/
$\kaz u\in\brn T{n}(\text{there is a unique\/
$v\in\brn S{n}$ with $u\sq v$})\,;$

\itlb{2sq4}%
if\/ $S\sq'_n T$ then\/ $S\sq_n T$ iff\/
$\kaz u\in\brn T{n-1}\,(\suc Tu=\suc Su)$.
\qed
\een
\ele

\ble
\lam{sost}
Let\/ $T\in\pew$, $n<\om$.
Assume that if\/ $u\in\brn T{n}$ then\/ $T_u\in\pew$,
$T_u\sq\req Tu$. 
Then the tree\/ $S=\bigcup_{u\in\brn T{n}}T_u$
belongs to\/ $\pew$ and satisfies\/ 
$S\sq_n T$ and $\req Su=T_u$ for all\/ $u\in\brn T{n}$.
\qed 
\ele

Note that under the conditions of the lemma, if 
$u\in\brn T{n}$ then $u\sq\roo{T_u}$, and 
in addition 
$\brn S{n}=\ens{\roo{T_u}}{u\in\brn T{n}}$.

\ble
\lam{fus}
Assume that\/
$\ldots \nq 4 T_4\nq 3 T_3\nq 2 T_2\nq 1 T_1\nq 0 T_0$ 
is an infinite decreasing sequence of trees in\/ $\pew$.
Then the tree\/ $T=\bigcap_nT_n$ belongs to\/ $\pew$, 
and we have\/ $T\sq_{n} T_n$, 
and hence\/ $\brn T{n}=\brn {T_{n+1}}{n}$, 
for all\/ $n$.\qed
\ele

\parf{Wide tree forcing notions and dense sets} 
\las{wds}

A non-empty set $\dP\sq\pew$ is a
\rit{wide tree forcing}, \sfo\ in brief,
\index{wide tree forcing}%
\index{forcing!wide tree forcing}%
\index{WT-forcing}%
\index{forcing!WT-forcing}%
if we have $\req Tu\in \dP$ whenever $u\in T\in\dP$. 
Thus $\pew$ itself is a \sfo, and if $S\in\pew$ then the set 
$
\ens{\req St}{t\in S}$ 
is a \sfo.

\bre
\lam{Pfo}
Any \sfo\ $\jpi$ can be considered as a forcing notion  
ordered so that if $T\sq T'$, then $T$ is a stronger
\usl. 
\index{forcing!stronger condition}%
The forcing $\jpi$ adjoins a cofinal element $x\in\on$. 
More exactly if a set $G\sq\jpi$ is 
\dd\jpi generic over a given set universe $\rV$ 
(and $\jpi\in\rV$ is assumed)  
then the intersection $\bigcap_{T\in G}[T]$ contains a  
unique element $\xx G \in (\omiv){}^\om$, 
and $\xx G$ satisfies  
\index{zzaG@$\xx G$, generic element}%
\index{generic element!$\xx G$}%
$G=\ens{T\in\jpi}{\xx G\in[T]}$, 
$\rV[G]=\rV[\xx G]$, 
and $\tsup{\xx G}=\omi^\rV$ (cardinality collapse). 

Elements $\xx G$ of this kind are called  
\dd\jpi\rit{generic}. 
\ere

To prove Theorem~\ref{Tun} we'll make use of
a certain \sfo\ $\dP\sq\pew$.

\bdf
\lam{denp}
A set\/ $D\sq{\jpi}$ is \rit{dense} in ${\jpi}$  
\index{set!dense}%
\index{dense}%
if for any $S\in{\jpi}$ there is a tree $T\in D$, 
$T\sq S$, 
\index{set!open dense}%
\index{open dense}%
\rit{open dense}, if in addition   
$S\in D$ holds whenever $S\in{\jpi}$, $T\in D$, 
$S\sq T$, and 
\index{set!pre-dense}%
\index{pre-dense}%
\rit{pre-dense}, if the set 
$D'=\ens{S\in\jpi}{\sus T\in D(S\sq T)}$ is dense.
\edf

\vyk{
\bdf
\lam{inext}
Assume that\/ $T\in\pew$, $D\sq\pew$, $t\in\bsd$, $H\sq\bsd$. 
Define
\bit
\item 
$T\ine D$ iff $T\sq S$ for some $S\in D$ --- 
let $\hdt DT=\ens{s\in T}{\req Ts\ine D}$;%
\imar{ine}%
\index{zzinext@$\ine$}%

\item
$t\ine H$ iff $s\sq t$ for some $s\in H$.\qed
\index{zzDIT@$\hdt DT$}%
\eit
\eDf
}

If $T\in\pew$ and $D\sq\pew$ then let
$\hdt DT=
\ens{s\in T}{\sus S\in D\,(\req Ts\sq S)}$.
\index{zzDIT@$\hdt DT$}%
\imar{hdt D T}%

\ble
\lam{plo1}
Assume that\/ $\jpi$ is a \sfo, and\/
$D_n\sq\jpi$ is pre-dense in\/ $\jpi$ for all\/ $n$.
Let\/ $S_0\in\jpi$.
Then there is a tree\/ $T\in\pew$
{\rm(not necessarily in $\jpi$!)}
such that\/ $T\sq S_0$ and if\/ $n<\om$
then\/ $\brn Tn\sq\hdt {D_n}T$.
\ele
\bpf
We wlog assume that each $D_n$ is open dense; 
otherwise replace it by 
$D'_n=\ens{S'\in\jpi}{\sus S\in D_n\,(S'\sq S)}$.
Using Lemma~\ref{sost} and the open density, 
define a sequence
$\ldots \nq 4 T_4\nq 3 T_3\nq 2 T_2\nq 1 T_1\nq 0
T_0\sq S_0$, such that if $n<\om$ and
$s\in\brn {T_{n+1}}{n}$ then $\req{T_{n+1}}s\in D_n$.
By Lemma~\ref{fus}, the tree $T=\bigcap_nT_n$
is as required: if $s\in\brn T{n}$ then 
$s\in\brn {T_{n+1}}{n}$, so that 
$\req{T}u\sq \req{T_{n+1}}u\in D_n$.
\epf

There is no way to directly extend Lemma~\ref{plo1} 
to the case of \dd\omi sequences of dense sets.
But a somewhat weaker result of Lemma~\ref{plo2} 
will be possible.

\parf{Bounded sets and continuous maps} 
\las{bscm}

It is known from descriptive set theory that if a
continuous map $f:P\to\bn$ is defined on a perfct
set $P\sq\bn$ then $f$ is a bijection or a constant
on a suitable perfect subset $P'\sq P$.
A similar but somewhat more complicated dichotomy holds
for wide trees.
Say that a set $X\sq\on$ is \rit{bounded},
if there is an ordinal 
\index{set!bounded}%
\index{bounded}%
$\ba<\omi$ such that $X\sq\ba^\om$. 
Note that if $T\in\pew$ then the set $[T]$ is
{\ubf un}bounded. 

\ble
\lam{minl}
Let\/ $S\in\pew$ and\/ $f:[S]\to\on$ be continuous. 
There is a tree\/ $T\sq S$, $T\in\pew$, such 
that either\/ $\ima f{[T]}$ is bounded 
or\/ $f\res[T]$ is a bijection.
\ele
\bpf
Suppose that for no $T\in\pew$, $T\sq S$, 
$f\res[T]$ is bounded. 
Then, as the set $\brn S1$ is uncountable, 
by a simple cardinality argument there exist: 
an uncountable set $U\sq\brn S1$,
a number $k$, 
and for each $t\in U$ --- an ordinal $\xi_t<\omi$ 
and a tree $U_t\in\pew$ satisfying 
$U_t\sq\req St$, 
$f(x)(k)=\xi_t$ for all $x\in[U_t]$
(same $k$ for all $t\in U$!), 
and if $t\ne t'$ belong to $U$ then 
$\xi_t\ne \xi_{t'}$.

Then the tree $S_1=\bigcup_{t\in U}U_t$ belongs
to $\pew$ and satisfies $S_1\sq'_1 S$.
In addition, there is a number $k=k_1$ such that
if $u\ne u'$ belong to $\brn{S_1}1$ and
$x,x'\in [S_1]$, $u\su x$, $u'\su x'$, then 
$f(x)(k_1)\ne f(x')(k_1)$.

Similarly, there is a tree $S_2\in\pew$,
$S_2\sq'_2 S_1$, and a number $k_2$, such that
if $u\ne u'$ belong to $\brn{S_2}2$ and
$x,x'\in [S_1]$, $u\su x$, $u'\su x'$, then 
$\ang{f(x)(k_1),f(x)(k_2)}\ne\ang{f(x')(k_1),f(x')(k_2)}$.

Iterating this construction appropriately by
induction, we get
a required
tree $T=\bigcap_n S_n\in\pew$ by Lemma~\ref{sost}.
\epf

The next theorem presents a dichotomy somewhat
different than the one considered by Lemma~\ref{minl},
and related to the
case of \dd\ali many maps.

\bdf
\lam{huf}
If $f:\on\to\on$ is a continuous map, 
and $U,V\in\pew$, then   
$\huf UfV$ is the set of all strings $s\in U$ 
\index{zzHufv@$\huf UfV$}%
such that
(1)
$[V]\cap(\ima f{[\req Us]})$ is bounded or
(2) 
$f\res[\req Us]$ is a \rit{total identity}, that is, 
\index{total identity}%
\index{identity}%
$f(x)=x$ for all $x\in [\req Us]$.

Note that (1) and (2) are incompatible provided
$U\sq V$.
\edf


\bte
[the proof ends in Section~\ref{pret}]
\lam{frt}
Assume that\/ $S\in\pew$ and, 
for each\/ $\al<\omi$, 
$f_\al:\on\to\on$ is a continuous function. 
Then there is a tree\/ $T\in\pew$, 
$T\sq_1 S$, such that, 
for any\/ $\al<\omi$, 
the set\/ $\huf T{f_\al}T$ 
is dense in\/ $T\,.$ 
\ete

\parf{Iteration of wide trees} 
\las{ite}

Here we develop another method of construction of 
trees in $\pew$, similar to a 
construction introduced in \cite[2.14]{abr2}, 
and designed for the proof of Theorem~\ref{frt}.

\bdf
\lam{dri}
A function $J$ is an \rit{iteration (of wide trees)},
\index{iteration!of wide trees}%
\index{iteration!core}%
\index{iteration!$\iwt$}%
\index{zziwt@$\iwt$}%
in symbol $J\in\iwt$, 
if its \rit{core} $\cn=\dom J$ is a subtree of $\bsd$  
(possibly with \eno s and/or isolated branches), 
all values $J(u)$ are trees 
in $\pew$, and in addition 
\ben
\nenu
\itlb{dri1}%
if $u\sq v$ belong to $C$ then $v\in J(u)$ 
and $J(v)\sq \req{J(u)}v$; 

\itlb{dri2}%
if $u\su v$ belong to $C$, 
$\lh v=\lh u+1$, and $u\nin\br{J(u)}$ 
then $J(v)=J(u)$. 

{\noii 
In this case we define the \rit{wrap} of $J$,
\index{iteration!wrap, $\tre J$}%
\index{zzwrJ@$\tre J$}%
}
\itlb{dri3}%
$\tre J=\ens{s\in\bsd}
{\kaz u\in\dom J\,(u\su s\imp s\in J(u)}$.  
\een
If $\jpi\sq\pew$ then $\ms\jpi$ consists 
of all iterations $J\in\iwt$ with $\ran J\sq\jpi$.
\index{iteration!$\ms\jpi$}%
\index{zziwtp@$\ms\jpi$}%
\imar{ms jpi}%
An iteration $J\in\iwt$ is \rit{small} if 
the core $C=\dom J$ is
at most countable. 
\index{iteration!small}%
\edf 

In particular $\pu\in\iwt$ and $\tre\pu=\bsd$.

If ${\cn}\sq\bsd$ is a tree and $s\in\bsd$ 
then let $\pro {\cn}s$ (the \rit{projection}) 
\imar{pro Cs}%
\index{zzprojCs@$\pro {\cn}s$}%
\index{projection!projCs@$\pro {\cn}s$}%
be the largest 
string in ${\cn}$ with $u\sq s$; $\pro {\cn}s=s$ 
provided $s\in {\cn}$. 

\ble
\lam{iwf}
If\/ $J\in\iwt$ 
then\/ $T=\tre J\in\pew$,  ${\cn}=\dom J\sq T$, 
and 
\ben
\renu
\itlb{iwf1}%
if\/ $s\in {\cn}$ then\/ $s\sq\roo{J(s)}$, 
$\req Ts\sq J(s)$, and\/ 
$\suc Ts=\suc{J(s)}s\,;$ 

\itlb{iwf2}%
if\/ $s\in T\bez {\cn}$ and\/ $u=\pro {\cn}s$ then\/ 
$s\in J(u)$ and\/ $\req Ts=\req{J(u)}s\;;$ 

\itlb{iwf3}%
if\/ $s\in {\cn}$ is an \eno\ in\/ ${\cn}$ then we have\/ 
$\req Ts=J(s)$.      
\een
\ele

\bpf
If $u\in {\cn}$ then $u\in T$ by \ref{dri}\ref{dri1}, 
so we have ${\cn}\sq T$. 

\ref{iwf1}
If $s\in {\cn}$ then $J(s)=\req{J(s)}s$ 
by \ref{dri}\ref{dri1} with $u=v=s$, so that 
obviously $s\sq\roo{J(s)}$. 
If now $t\in T$ and $s\sq t$ then $t\in J(s)$ 
by \ref{dri}\ref{dri3}, therefore 
$\req Ts\sq J(s)$. 
This implies $\suc Ts\sq\suc{J(s)}s$. 
To get the equality, let 
$t=s\we \xi\in \suc{J(s)}s$. 
Then $t\in T$ by \ref{dri}\ref{dri1},\ref{dri3},
so $t\in \suc{T}s$, 
as required.

\ref{iwf2}
If $s\nin {\cn}$ then by \ref{dri}\ref{dri1}\ref{dri3} the 
criterion of $s\in T=\tre J$ is just $s\in J(u)$, where 
$u=\pro {\cn}s$. This easily implies the result. 
And 
\ref{iwf3}
is similar to \ref{iwf2}.

To prove $T\in\pew$, let $s\in T$. 
We have to prove that 
(a) if $s\in \br T$ then $\suc Ts$ is 
uncountable, and 
(b) there is a string $s'\in \br T$ with  
$s\sq s'$.
By \ref{iwf1}, \ref{iwf2} we have (a) immediately, 
so it remains to check (b).\vom

{\sl Case 1\/}:
$s\in T\bez {\cn}$. 
Then 
$\req Ts=\req{J(u)}s$ by \ref{iwf2}, where $u=\pro {\cn}s$. 
But $\req{J(u)}s\in\pew$ by Lemma~\ref{iso}, which 
easily implies (b).\vom
 
{\sl Case 2\/}:
$s$ is an \eno\ in ${\cn}$, 
so $\req Ts=J(s)\in\pew$ by \ref{iwf3}, 
follow Case 1.\vom

{\sl Case 3\/}:
there is an \eno\ $s'$ in ${\cn}$ with $s\sq s'$ ---
apply Case 2 for $s'$.\vom

{\sl Case 4\/}:  
if all the above fails then there is an infinite
branch in ${\cn}$ containing $s$, that
is, $b\in\on$ such that $b\res n\in {\cn}$, $\kaz n$,
and $s=b\res n_0$, where $n_0=\lh s$.
Then $b\res n\in J(s)$ for all $n$ by \ref{iwf1}.
Therefore, as $J(s)\in\pew$, there is a least number
$k\ge n_0$ with $t=b\res k\in\br{J(s)}$.
Then by the way $J(t)=J(s)$ by \ref{dri}\ref{dri2},
hence $t\in\br{J(t)}$, and finally $t\in\br{T}$
by \ref{iwf1}, as required.
\epf

\vyk{
{\sl Case 3\/}:
$s\in {\cn}$ is not an \eno\ of ${\cn}$. 
If there is an \eno\ $s'$ in ${\cn}$ with $s\su s'$ then 
(b) follows from the result in Case 2. 
If there is a node  $s'\in \br {\cn}$ with $s\su s'$  
then $s'\in \br T$ as well, so we have (b).
The following is the remaining case in the proof 
of (b).\vom 

{\sl Case 4\/}:  
the set ${\cn}_s=\ens{u\in {\cn}}{s\sq u}$ 
is a single branch $b\in\on$, thus $s=b\res n_0$, 
where $n_0=\lh s$, 
and if $n>n_0$ then $b\res n\in {\cn}$.\vom

{\sl Subcase 4a\/}: 
there is $n\ge n_0$ such that 
$b\res n\in\br{J(b\res n)}$.
Then $b\res n\in\br T$ by \ref{iwf1}, and 
$s\sq b\res n$, as required.\vom

{\sl Subcase 4b\/}: 
$b\res n\nin\br{J(b\res n)}$ for all $n\ge n_0$.
Then $J(b\res n)=J(s)$ for all $n\ge n_0$, 
by Definition~\ref{dri}\ref{dri2}. 
It follows that $b$ is an infinite branch in 
$J(s)\in\pew$. 
Therefore there exists a number $n\ge n_0$ 
satisfying $b\res n\in\br{J(s)}=\br{J(b\res n)}$.
We conclude that $b\res n\in\br{T}$ by \ref{iwf1}, 
as required.
} 

The lemma allows to maintain infinite, even 
uncountable \dd\sq decreasing   
sequences of trees in $\pew$, 
with the help of the following two rather obvious results.

\ble
\lam{itmo}
If\/ $J\sq J'$ are  iterations in\/ $\iwt$ then\/ 
$\tre{J'}\sq \tre{J}$. 
\lam{iteu}

If\/ $\sis{J_\xi}{\xi<\la}$ is a\/ \dd\su increasing 
sequence of iterations\/ $J_\xi\in\iwt$ then\/ 
$J=\bigcup_{\xi<\la}J_\xi\in\iwt$ and\/ 
$\tre{J}=\bigcap_\xi\tre{J_\xi}$.\qed
\ele

\ble
\lam{+1}
Let\/ $J\in\iwt$, $\dom J={\cn}\sq {\cn}'\sq T=\wra{J}$,
${\cn}'$ be a tree. 
\bit
\item
{\rm 
Define a\/ {\ubf natural extension} 
$J'
$ of $J$ 
%
to 
${\cn}'$ 
by $\dom{J'}={\cn}'$,
$J'(s)=J(s)$ for $s\in {\cn}$, and if\/ 
$s\in {\cn}'\bez {\cn}$ and $u=\pro {\cn}s$ 
then $J'(s)=\req{J(u)}s$.}
\eit 
Then\/ $J'\in\iwt$, $J\sq J'$, 
$\tre {J'}=\tre J$.\qed
\ele

Condition \ref{dri2} of Definition~\ref{dri} imposes
important restrictions on the construction of
iterations, basically justifying proper shrink only
at successors of branching nodes.
Nevertheless it leaves us enough freedom.

\ble
\lam{+2}
Assume that\/ $\jpi$ is a \sfo, $J\in\ms\jpi$, 
${\cn}=\dom J$, 
$s\in T=\tre J$, and\/ $s\nin {\cn}$ or\/ 
$s$ is an \eno\ in\/ ${\cn}$.
Let\/ $U\in\jpi$, $U\sq \req Ts$. 
Then there exists an iteration\/ $J'\in\ms\jpi$ 
and a string\/ $s'\in\dom J'$ 
such that\/ $J\sq J'$, $s\sq s'$, 
and\/ $J'(s')\sq U$.
\ele
\bpf
Let $t=\roo U$, thus $s\sq t\in\br U$ and 
all shorter strings $v\su t$ do not 
belong to $\br U$.
Pick any $s'\in U$ with $\lh{s'}=\lh t+1$; 
then $t\su s'\nin {\cn}$.
Let $u=\pro {\cn}s$.
Let 
$J'\in\ms\jpi$ be the extension of $J$ to the 
domain 
${\cn}'={\cn}\cup\ens{v}{u\su v\sq t}\cup\ans{s'}$ by 
$J'(u)=J(s)=\req{J(s)}u$ whenever $s\su u\sq t$, 
and finally $J'(s')=\req{U}{s'}$.
To see that \ref{dri}\ref{dri2} is satisfied 
for $J'$ at $u=t$ and $v=s'$, recall that 
$t\in \br U$, hence $t\in\br{J(s)}=\br{J(t)}$ 
as well.
\epf


\parf{Key dichotomy lemma} 
\las{nef}

\ble
\lam{rd}
Assume that\/ $\jpi\sq\pew$ is a \sfo, 
$J\in\ms\jpi$ is a small iteration, 
$S=\tre J$,   
$g_0\in {\cn}=\dom J$, 
and\/ $f:\on
\to\on$ is continuous. 
There is a small iteration\/ $J'\in\ms\jpi$ 
and a string\/ 
$g\in {\cn}'=\dom {J'}$, such that\/ $g_0\sq g$, 
$J\sq J'$, and\/ $g\in\huf TfT$, where\/ 
$T=\tre{J'}$, \ie,
\bce 
$(1)$  
$[T]\cap (\ima f[\req {T}{g}])$ is bounded, \quad 
{\ubf or} \quad $(2)$
$f$ is a total identity on\/ $[\req T{g}]$. 
\vyk{     
\item[or  \rm(3)] 
$[T]\cap(\ima f[\req {T}{g}])=\pu$. 
}%
\ece
\ele

The lemma will be crucial 
in the proof of Theorem~\ref{frt} in Section~\ref{pret}.

\bpf
Pick any $g_1\in S\bez {\cn}$ satisfying $g_0\sq g_1$. 
Let $u=\pro{{\cn}}{g_1}$. 
If $\ima f[\req {S}{g_1}]\sq[{\cn}]$ (a bounded set) then  
let $J_1\in\ms\jpi$ be the natural 
extension of $J$ to the domain
${\cn}_1={\cn}\cup\ens{s}{u\su s\sq g_1}$ by 
Lemma~\ref{+1}. 
Thus $J\sq J_1$, 
$\dom{J_1}={\cn}_1$, $J_1(g_1)=\req{J(u)}{g_1}$, and  
$\wra{J_1}=S$.
Therefore $J'=J_1$ and $g=g_1$ satisfy (1). 

Thus suppose that $x_1\in[\req {S}{g_1}]$, and 
$y_1=f(x_1)\in [S]\bez[{\cn}]$. 
As $f$ is continuous while $\on\bez [{\cn}]$ open, 
there is a longer string $g_2\in S\bez {\cn}$, 
$g_1\su g_2$, such that $f(x)\nin [{\cn}]$ for all 
$x\in\req S{g_2}$. 
If $f$ is a total identity on\/ $[\req T{g_2}]$ 
then let $J_2\in\ms\jpi$ be the natural 
extension of $J$ to the domain
${\cn}_2={\cn}\cup\ens{s}{u\su s\sq g_2}$ by 
Lemma~\ref{+1}; now $J'=J_2$ and $g=g_2$ satisfy (2). 

Thus suppose that $x_2\in[\req {S}{g_2}]$, and 
$y_2=f(x_2)\ne x_2$. 
There is a yet longer string 
$g_3\in S\bez {\cn}$, 
$g_2\su g_3$, such that $f(x)\ne x$ 
and $f(x)\nin [{\cn}]$ for all 
$x\in\req S{g_3}$. 
If $(\ima f{[\req S{g_3}]})\cap[S]=\pu$ 
then let $J_3\in\ms\jpi$ be the natural 
extension of $J$ to the domain
${\cn}_3={\cn}\cup\ens{s}{u\su s\sq g_3}$;  
now $J'=J_3$, $g=g_3$ satisfy (1). 

Thus suppose that $x_3\in\req S{g_3}$ and 
$y_3=f(g_3)\in[S]$. 
In addition, $x_3\ne y_3\nin [{\cn}]$ holds as 
$g_2\sq g_3$, 
hence there is $m\ge\lh{g_3}$ such that 
$t=y_3\res m\in\br S\bez {\cn}$ 
and $t\ne s=x_3\res m$. 
Let $t'= y_3\res{(m+1)}$ 
(a successor of $t$ in $S$). 
There is a string $h\in S$ 
such that $t'\su h$ but $h\ne t''=y_3\res\ell$, 
where $\ell={\lh h}$. 
As $f$ is continuous, pick a number 
$n\ge n_3=\lh{g_3}$ 
such that $t''\su f(x)$ 
holds for all $x\in[\req S{g}]$, 
where $g=x_3\res n$. 
Recall that $u=\pro {\cn}g$. 
Let $v=\pro {\cn}t$,\pagebreak[0] 
$$
{\cn}'={\cn}
\cup\ens{w\in\bsd}{u\su w\sq g} 
\cup\ens{w\in\bsd}{v\su w\sq t'}\,,
$$ 
and extend the iteration $J$ to the domain ${\cn}'$ by 
$J'(w)=\req{J(u)}w$ whenever $u\su w\sq g$,  
$J'(w)=\req{J(v)}w$ whenever $v\su w\su t$, and 
$J'(t')=\req{J(v)}{h}$. 

Now it suffices to prove (1) in the form  
$[T]\cap(\ima f[\req {T}{g}])=\pu$. 
Let $g\su x\in[S]$.  
Then $y=f(x)$ satisfies $t''\su y$, 
hence $h\not\su y$. 
Let's show that $y\nin[T]$. 
It suffices to check $t''\nin T$. 
Suppose otherwise. 
Then, as $t'\in {\cn}'$, we have 
$t''\in J'(t')$ by \ref{dri}\ref{dri3}.
However $J'(t')= \req{J(v)}{h}$, 
so it follows that $t''$ and $h$ are compatible, 
which contradicts to the construction, as required.
\epf

\parf{The proof of the restriction theorem} 
\las{pret}

Here we accomplish the proof of Theorem~\ref{frt} 
on the base of the results above. 
We argue in the assumptions of Theorem~\ref{frt}. 

The set $\jpi=\ens{\req St}{t\in S}$ is a \sfo\ 
and $S\in\jpi$.

By Lemmas \ref{rd} and \ref{itmo}, \ref{+1}, 
there is a  
\dd\sq increasing sequence of small iterations 
$J_\ga\in\ms\jpi$, $\ga<\omi$, with domains 
${\cn}_\ga=\dom{J_\ga}$ and trees $S_\ga=\wra{J_\ga}$, 
such that ${\cn}_0=\ens{u}{u\sq \sg}$, where $\sg=\roo S$, 
and $J_0(u)=S$ for all $u\in {\cn}_0$, 
the sets ${\cn}=\bigcup_{\ga<\omi}{\cn}_\ga$ and 
$T=\bigcap_{\ga<\omi}T_\ga$ coincide 
(Lemma~\ref{+1} is responsible), 
and in addition 
(Lemma~\ref{rd} is responsible), 
if $s_0\in {\cn}=T$ and $\al<\omi$ then 
there is an index $\ga=\ga(s_0,\al)<\omi$ and a 
string $s\in {\cn}_\ga$ such that $s_0\sq s$ and 
$s\in\huf{T_\ga}{f_\al}{T_\ga}$.
%
Then $J=\bigcup_\al J_\al\in\ms\jpi$,
${\cn}=\dom J$, and $T=\tre J$, by Lemma~\ref{iteu}.
Moreover, as $T\sq T_\ga$, 
we have $s\in\huf{T}{f_\al}{T}$ as well. 
It follows that  
the set 
$\huf{T}{f_\al}{T}$
is dense in $T$, and obviously \rit{open} dense.
And finally we have $T\nq 1 S$ by
Lemma~\ref{iwf}\ref{iwf1}
with $s=\sg=\roo S$. (Recall that $J_0(\sg)=S$.) 

\vyk{
Prove that $T\nq 1 S$, not merely $\sq$. 
It suffices to check that if $\ga<\omi$ and 
$\sg\we\ga\in S$ then $\sg\we\ga\in T$,
where $\sg=\roo S$. 
Indeed if $\sg\we\ga\in {\cn}$ then $\sg\we\ga\in T$ as 
${\cn}\sq T$. 
If $\sg\we\ga\nin {\cn}$ then the only reason for 
$\sg\we\ga\nin T$ is $\sg\we\ga\nin J_\xi(\sg)$ for 
some $\xi$, 
but we have $J_\xi(\sg)=J_0(\sg)=S$ for all $\xi$ 
by construction.
}

\vom

\qeD{Theorem~\ref{frt}}

\parf{Belts and covering} 
\las{bc}

Here we introduce the last major tool employed in the 
definition of the forcing notion for Theorem~\ref{Tun}. 
It is based on the following definition. 

\bdf                               
\lam{mbelt}
A set $H\sq\bsd$:
\bit
\item[$-$]
\rit{meets} $x\in\on$ iff 
\index{meets}
$\sus m\,(x\res m\in H)$; 
\item[$-$]
is a \rit{belt} for 
a tree 
$T\in\pew$, 
if it meets every $x\in[T]$;

\item[$-$]
\rit{weakly covers} $T$, 
in symbol $T\sqw B$, 
if there is an ordinal $\ba<\omi$ 
such that $H$ is a belt for each subtree $\req Ts$, where 
$s\in T$ and $\tmax s\ge\ba$ --- in other words, we require 
$H$ to meet every $x\in[T]$ with  $\tsup x\ge\ba$.
\eit
\index{belt}%
\index{weakly covers}%
\index{zzTwD@$T\sqw D$}%
For instance, if $n<\om$ then $\brn Tn$ is a belt for $T\in\pew$.
\edf

\ble
\lam{lwec}
Let\/ $H\sq T$ weakly cover\/ $T\in\pew$  
with a parameter $\ba<\omi$.
Then 
\ben
\renu
\itlb{lwec0}%
$H$ is pre-dense in\/ $T\,;$

\itlb{lwec1}%
$H$ weakly covers any tree\/ $S\in\pew$, $S\sq T$, 
with the same\/ $\ba\,;$

\itlb{lwec2}%
the set\/ 
$ 
X=\ens{x\in[T]}{H\text{\rm\ does not meet }x}
$  
satisfies\/ $X\sq\ba^\om$.
\een 
\ele
\bpf
\ref{lwec2}
Let $x\in [T]\bez\ba^\om$, $x(j)\ge\ba$ 
for some $j$. 
Let $s=x\res{(j+1)}$. 
Then $H$ is a belt for $\req Ts$, hence 
$H$ meets $x$.
\epf

\bre
\lam{belta}
Being a belt  
is equivalent to the wellfoundedness of the subtree 
$T'=\ens{s\in T}{\neg\:\sus t\in H\,(t\sq s)}$, 
hence it is an \rit{absolute} notion. 
It follows that
to weakly cover with a parameter $\ba$ 
is an \rit{absolute} notion, too. 

Now assume that $H\sq T$ weakly covers $T\in\pew$ with 
a parameter $\ba<\omi$. 
Let $x\in[T]$ be an element cofinal in $\omi$ 
(=$\omiv$ of the given set universe $\rV$), 
which may exist in an extension of $\rV$, 
Remark~\ref{Pfo}. 
{\ubf We claim that\/ $H$ meets\/ $x$.} 
Indeed, $x\nin\ba^\om$ by the cofinality, 
and on the other hand, the absoluteness of the weak 
covering allows 
to apply Lemma~\ref{lwec}\ref{lwec2} in the extension 
containing $x$.  
\ere

\ble
\lam{plo2}
Assume that\/ $\jpi$ is a \sfo, $T\in\jpi$, and\/
$D_\xi\sq\jpi$ is open dense in\/ $\jpi$ for all\/
$\xi<\omi$.
Then there is a tree\/ $S\in\pew$ such that\/ 
$S\nq1 T$,
and each set\/ 
$\hdt{D_\xi}{S}=
\ens{t\in S}{\sus U\in D_\xi\,(\req St\sq U)}$ 
weakly covers\/ $S$.
\ele

\vyk{Recall: 
$\hdt DT=\ens{s\in T}{\sus S\in D\,(\req Ts\sq S)}$  
whenever
$T\in\pew$, $D\sq\pew$.
\index{zzDIT@$\hdt DT$}%
}

\bpf
If $\al<\omi$ then fix an enumeration of the
countable set
$\ens{D_\xi}{\xi\le\al}=\ens{D^\al_k}{k<\om}$.
Using Lemma~\ref{sost} and the open-density of each
$D_\xi$ in $\jpi$, define a sequence
$\ldots\nq 5 T_4\nq 4 T_3\nq 3 T_2\nq 2 T_1\nq 1 T_0=T$ 
of trees in\/ $\pew$, such that if $n\ge1$ and
$u\in\brn{T_{n}}{n}$ then
$\req{T_{n}}u\in \bigcap_{j,k\le n} D^{u(j)}_{k}$.
The tree $S=\bigcap_nT_n$ belongs to $\pew$ and 
satisfies $S\sq_{n+1} T_n$ and  
$\brn {S}n=\brn {T_{n}}n$ for all $n$, 
by Lemma~\ref{fus}. 
In particular $S\nq1 T$.
Now suppose that $\xi<\omi$.

We claim that $\xi$ itself witnesses $\hdt{D_\xi}{S}$ 
to weakly cover $S$.
Let $x\in[S]$ and $x(j)=\al\ge\xi$ for some $j$. 
Then $D_\xi=D^{\al}_k$ for some $k$. 
Let $n=1+\max\ans{j,k}$. 
There is a number $m\ge n$ such that $u=x\res m$ belongs 
to $\brn {S}n=\brn {T_{n}}n$. 
Then 
$\req{S}u\sq\req{T_{n}}u\in D^{u(j)}_{k}=D^{\al}_{k}=D_\xi$
by construction, and 
we are done.
\epf 

\bcor
\lam{plo3}
If\/ $T\in\pew$ and\/
$H_\xi\sq T$ is open dense in\/ $T$ for all\/
$\xi<\omi$ then there is\/ $S\in\pew$ such that\/ 
$S\nq1 T$ and each\/ $H_\xi\cap S$ 
weakly covers\/ $S.$ 
\ecor
\bpf
Apply the lemma for $\jpi=\ens{\req Ts}{s\in T}$ and 
$D_\xi=\ens{\req Ts}{s\in H_\xi}$.
\epf

\parf{Extensions of wide tree forcing notions} 
\las{wtfn}

The forcing notion to prove Theorem~\ref{Tun}
will be defined in the form of an \dd\omi union 
of its parts --- \sfo s of cardinality $\le\ali$.

\bdf
\lam{dfex}
Let $\cM$ be any set and $\jpi$ be a \sfo.  
Another \sfo\ $\jqo$ 
is an \dd\cM\rit{extension} of $\jpi$, 
in symbol $\jpi\prol\cM\jqo$, if the following holds:%
\index{forcing!\dd\cM extension, $\prol\cM$}%
\index{Mextension@\dd\cM extension, $\prol\cM$}%
\index{zzsqsubM@$\prol\cM$}%
\index{zzsqsub@$\prok$}%
\imar{prol cM}%
\ben
\Aenu
\itlb{dfex1}%
$\jqo$ is dense in $\jqo\cup\jpi$; 

\itlb{dfex2x}%
$\jqo$ \rit{refines} $\jpi$: 
if $Q\in\jqo$ then there exists $T\in\jpi$ 
satisfying $Q\sq T$;

\itlb{dfex3}%
if a set $D\in\cM$, $D\sq\jpi$ is pre-dense in $\jpi$ 
and $U\in\jqo$ then the set 
$\hdt DU=\ens{s\in U}{\sus S\in D\,(\req Us\sq S)}$
weakly covers $U$;

\itlb{dfex4}%
if $T_0\in\jpi$ and
$\sis{D_n}{n<\om}\in\cM$ is a sequence of pre-dense 
sets $D_n\sq\jpi$ then there is a 
tree $T\in\jqo$ such that $T\sq T_0$, and
$\brn Tn\sq \hdt{D_n}T$ for all 
$n$;

\itlb{dfex5}%
if $T_0\in\jpi$ and $f:\on\to\on$, $f\in\cM$, 
is continuous, then there is $T\in\jqo$ 
such that $T\sq T_0$, and 
either $\ima f{[T]}$ is bounded or
$f\res [T]$ is a bijection;

\itlb{dfex6}%
if $f\in \cM$, $f:\on\to\on$ is a continuous map, 
and $U,V\in\jqo$, 
then 
the set $\huf UfV$, of all strings $s\in U$  
such that $[V]\cap(\ima f{[\req Us]})$ is bounded or 
$f\res[\req Us]$ is a total identity, weakly covers 
$U$. 
%
\qed 
\een
\eDf

If $\cM=\pu$ then we write $\jpi\prok\jqo$ 
\index{forcing!extension, $\prok$}%
\imar{prok}%
instead of $\jpi\prol\pu\jqo$; 
in this case \ref{dfex3} -- \ref{dfex6}
are trivial.
Generally, in the role of $\cM$ we'll consider  
transitive models of the theory $\zfc'$ 
\index{zzzfcp@$\zfc'$}%
\index{theory!zfcp@$\zfc'$}%
\index{model!\dd\sg closed}%
which includes all $\ZFC$ axioms except for the Power   
Set axiom, but an axiom is adjoined, which claims  
the existence of $\omi$ and $\pws\omi$. 
(Then the existence of sets like $\bsd$ and 
$\pew$ easily follows.)

\vyk{
We'll mostly consider \dd\sg\rit{closed} 
transitive models 
$\cM\mo\zfc'$, \ie, those satisfying: if 
$f:\om\to\cM$ then $f\in\cM$. 
If $\cM$ is such then $\omi\in\cM$.
}

\ble
\lam{teb}
Let\/ $\jpi,\jqo,\jR$ be\/ 
\sfo s satisfying\/ $\jpi\prok\jqo\land\jqo\prok\jR$. 
Then\/ $\jpi\prok\jR$, and if\/ {\rm{(K)}} is one of\/ 
{\rm\ref{dfex3}, \ref{dfex4}, \ref{dfex5}, \ref{dfex6}} 
and the pair\/ $\jpi\prok\jqo$ satisfies\/ {\rm{(K)}} 
with some\/ $\cM$, then the pair\/ $\jpi\prok\jR$ 
satisfies\/ {\rm{(K)}} with the same\/ $\cM$.
\ele
\bpf
\ref{dfex3} 
Let $R\in\jR$. 
As $\jqo\prok\jR$, there is a tree $Q\in\jqo$ with $R\sq Q$. 
Then $\hdt DQ$ weakly covers $Q$ by \ref{dfex3} for 
$\jpi,\jqo$. 
Then easily $\hdt DR$ weakly covers $R$. 


\ref{dfex4} 
If $T'\sq T$ and $t\in\brn{T'}n$ then there is 
a string $s\in\brn Tn$ with $s\sq s'$. 

\ref{dfex6}
If $U'\sq U$ and $V'\sq V$ then 
$\huf UfV\cap U'\sq \huf {U'}f{V'}$.
\epf

\ble
\lam{tec}
Assume that\/ $\cM\models\zfc'$ is a transitive model,    
and\/ $\jpi\in \cM$ and\/ $\jqo$ are\/ 
\sfo s satisfying\/ $\jpi\prol\cM\jqo$. 
Then
\ben
\renu
\itlb{tec1}%
if a set\/ $D\in\cM$, $D\sq\jpi$ 
is pre-dense in\/ $\jpi$ 
then\/ $D$ 
is pre-dense in\/ $\jpi\cup\jqo\;;$

\itlb{tec2}%
if\/ 
$T,T'\in\jpi$  are incompatible in\/ $\jpi$  
then\/ $T,T'$ are incompatible in\/
$\jpi\cup\jqo\;.$ 
\een
\ele
\bpf
\ref{tec1}
Let $U\in\jqo$. 
Then $\hdt DU$ weakly covers $U$ 
by \ref{dfex}\ref{dfex3}. 
Let $s\in\hdt DU$. 
Then $U'=\req Us\in\jqo$, $U'\sq U$, and $U'\sq S$ 
for some $S\in D$. 

\ref{tec2}
The sets 
$D(T)=\ens{S\in\jpi}{S\sq T\lor [S]\cap[T]=\pu}$
and $D(T')$ 
belong to $\cM$ and are open dense in $\jpi$ by 
Lemma~\ref{iso}.
Therefore $D=D(T)\cap D(T')$ is open dense either, 
and in fact 
$S\in D\imp [S]\cap[T]=\pu\lor [S]\cap[T']=\pu$ 
by the incompatibility.
It follows that if $U\in\jqo$ and, by \ref{tec1},  
$S\in D$ and $U'\in\jqo$, $U'\sq U$, 
$U'\sq S\cap U$,
then $[U']\cap[T]=\pu$ or $[U']\cap[T']=\pu$, 
hence $U$ cannot witness the compatibility of $T,T'$.
\epf 

We now establish the existence of extensions.

\bte
\lam{tex}
Assume that\/ $\cM\models\zfc'$ is a\/  
transitive model 
of cardinality\/ ${\le}\,\ali$,   
and\/ $\jpi\in \cM$ is a \sfo, 
$\card\jpi\le\ali$ in\/ $\cM$.  
Then there exists a\/ \sfo\/ $\jqo$
of cardinality\/ $\ali$, 
satisfying\/ $\jpi\prol\cM\jqo$.
\ete

\bpf
{\it Step 1}. 
If $P\in \jpi$ then by Lemma~\ref{plo2} there is 
a tree $T(P)\in\pew$, $T(P)\sq P$, such that 
$\hdt D{T(P)}$ weakly covers $T(P)$ for each 
$D\in\cM$, 
$D\sq \jpi$, predense in $\jpi$. 
The set 
$\jpi'=\ens{\req{T(P)}s}{P\in\jpi\land s\in T(P)}$  
is a \sfo\ of cardinality $\ali$ and  
\ref{dfex}\ref{dfex1},\ref{dfex2x},\ref{dfex3} 
hold for $\jqo=\jpi'$. 

{\it Step 2}.  
To fulfill \ref{dfex}\ref{dfex5}, if $P'\in\jpi'$   
and $f:\on\to\on$, $f\in\cM$ is continuous, then 
by Lemma~\ref{minl} there is a tree $T(P',f)\in\pew$, 
such that $T(P',f)\sq T$, and 
either $\ima f{[T(P',f)]}$ is bounded or
$f\res [T(P',f)]$ is a bijection. 
We let 
$\jpi''=\ens{\req{T(P',f)}s}{P'\in\jpi'\land s\in T(P',f)}$. 
Now \ref{dfex}\ref{dfex1},\ref{dfex2x},\ref{dfex3},\ref{dfex5} 
hold for $\jqo=\jpi''$. 
 
{\it Step 3}. 
To fulfill \ref{dfex}\ref{dfex4}, note first of all that 
each set $D\in\cM$, $D\sq\jpi$, pre-dense in\/ $\jpi$, 
remains pre-dense in\/ $\jpi\cup\jpi''$ by 
Lemma~\ref{tec}\ref{tec1}. 
If $P''\in\jpi''$ and $\jjd=\sis{D_n}{n<\om}\in\cM$ 
is a sequence of pre-dense sets $D_n\sq\jpi$, 
then by Lemma~\ref{plo1} there is a 
tree $T(P'',\jjd)\in\pew$ such that $T(P'',\jjd)\sq P''$, 
and if $n<\om$ and $s\in\brn {T(P'',\jjd)}n$ then
$\sus S\in D_n\,(\req {T(P'',\jjd)}s\sq S)$. 
We let 
$$
\jpi'''=\ens{\req{T(P'',\jjd)}s}
{P''\in\jpi''\land \jjd\in \cM}\,.
$$ 
Now 
\ref{dfex}\ref{dfex1},\ref{dfex2x},%
\ref{dfex3},\ref{dfex4},\ref{dfex5} 
hold for $\jqo=\jpi'''$. 
 
To fulfill \ref{dfex}\ref{dfex6}, we begin with some 
notation.
If $S\in\pew$ and $\al<\omi$ then let 
$\al\we S=\ens{\al\we s}{s\in S}$; 
then $\al\we S\in\pew$ and 
$\ang\al\sq\roo{\al\we S}$.
Conversely, if $W\in\pew$ and 
$\ang\al\sq\roo{W}$ 
then let $\niz {W}=\ens{s\in\bsd}{\al\we s\in W}$; 
then $\niz {W}\in\pew$ and $W=\al\we{(\niz W)}$.
We have $\niz{(\al\we S)}=S$, of course.

{\it Step 4}.  
Let 
$\jpi'''=\ens{R_\al}{\al<\omi}$. 
We convert $\jpi'''$ into a single tree 
$$
\textstyle
R=
\ans\La\:\cup\:\bigcup_{\al<\omi}(\al\we R_\al) 
\in\pew\,;
\quad  
\text{ then }\, 
R_\al=\niz{(\req R{\ang\al})}\,,\;\kaz \al.
$$ 
If $f:\on\to\on$ is continuous and $\al,\ba<\omi$ then 
define $f_{\al\ba}:\on\to\on$ so that  
$f_{\al\ba}(\al\we x)=\ba\we f(x)$, 
$f_{\al\ba}(\ba\we x)=\al\we f(x)$, 
and $f_{\al\ba}(y)=y$ whenever $y(0)\ne\al,\ba$. 
The set of continuous functions 
$F=\ens{f_{\al\ba}}{f\in\cM\land \al,\ba<\omi}$ 
is still of cardinality $\ali$. 
By Theorem~\ref{frt} there exists a tree 
$T\in\pew$, $T\nq1 R$, such that if $h\in F$ then 
the set $\huf ThT$ 
is open dense in $T$. 
Therefore by Corollary~\ref{plo3} there is a tree 
$Q\in\pew$ such that $Q\nq1 T$ 
(hence $Q\nq1 R$ as well) 
and if $h\in F$ then $\huf ThT$ weakly covers $Q$. 
Then $\huf QhQ$ weakly covers $Q$ as well by 
Lemma~\ref{lwec}\ref{lwec1} 
since  
$\huf ThT\cap Q\sq \huf QhQ$.

{\it Step 5}.  
Note that if $\al<\omi$ then 
the one-term string $\ang\al$ belongs to $Q$ since 
$Q\nq1 R$. 
Now let 
$Q_\al=\niz{(\req Q{\ang{\al}})}=
\ens{q\in\bsd}{\al\we q\in Q}$.
We claim that 
{\it
the \sfo\/ 
$\jqo=\ens{\req{Q_\al}q}{\al<\omi\land q\in Q_\al}$ 
satisfies\/ $\jpi\prol\cM\jqo$.}

First of all, 
$\jpi\prok{\jpi'}\prok{\jpi''}\prok{\jpi'''}\prok\jqo$
by construction, and hence $\jpi\prok\jqo$ holds, 
and we have 
\ref{dfex}\ref{dfex3},\ref{dfex4},\ref{dfex5} 
for the pair $\jpi\prok\jqo$ by Lemma~\ref{teb}.
\vyk{
\ref{dfex}\ref{dfex1} \
If $T\in\jpi$ then there are trees 
$T'\in\jpi'$, $T''\in\jpi''$, $T'''\in\jpi'''$   
with $T'''\sq T''\sq T'\sq T$. 
Then $T'''=R_\al$ for some $\al$, 
hence $Q_\al\in\jqo$ and $Q_\al\sq T$.

\ref{dfex}\ref{dfex3} \
We have 
$Q_\al\sq \niz{(\req R{\ang\al})}=R_\al\in\jpi'''$. 
By construction $R_\al\sq T''\sq T'$ holds, where 
$T'\in\jpi'$, $T''\in\jpi''$. 
Thus if a set $D\in\cM$ is pre-dense in $\jpi$ then 
$T'\sqw \hdt D{T'}=
\ens{s\in T'}{\req{T'}s\ine D}$ by Step~1, 
and hence $Q_\al\sqw \hdt D{Q_\al}$.

\ref{dfex}\ref{dfex4}\ref{dfex5} --- similar to the above.
}%


To check 
\ref{dfex}\ref{dfex6}, 
let $f\in \cM$, $f:\on\to\on$ be continuous, 
and $U=Q_\al$, $V=Q_\ba$ 
be trees in $\jqo$. 
To prove that $\huf UfV$ weakly covers $U$,    
let $h=f_{\al\ba}$. 
Then $\huf QhQ$ weakly covers $Q$ by Step~4. 
Thus there is an ordinal $\xi<\omi$ 
such that if $x\in [Q]$ and $\tsup x\ge\xi$ then 
$\huf QhQ$ meets $x$, so 
$x\res m\in  \huf QhQ$ for some $m$. 
We claim that $\xi$ witnesses that $\huf UfV$ 
weakly covers $U$. 

Assume that $y\in [U]=[Q_\al]$, $\tmax y\ge\xi$. 
Then $x=\al\we y\in [Q]$, 
so $s=x\res {(m+1)}\in \huf QhQ$ 
for some $m$, by the above.
Then $s=\al\we t$, where $t=y\res m$. 
It remains to prove that\/ $t\in \huf UfV$.

{\it Case 1\/}: 
$[Q]\cap (\ima h{[\req Qs]})$ is bounded. 
However $h=f_{\al\ba}$ and $U=Q_\al$, $V=Q_\ba$, hence 
$[Q]\cap (\ima h{[\req Qs]})=
\ba\we ([V]\cap(\ima f{[\req Ut]}))$. 
Thus the set $[V]\cap(\ima f{[\req Ut]})$ is bounded, 
therefore $t\in \huf UfV$.

{\it Case 2\/}: 
$h\res {[\req Qs]}$ is a total identity,   
$h(x)=x$ whenever $x\in \req Qs$. 
Then $\ba=\al$, $U=V$, and 
$f\res {[\req Ut]}$ is a total identity,    
thus still $t\in \huf UfV$.
\epf

\vyk{
To check 
\ref{dfex}\ref{dfex6}, 
let $f\in \cM$, $f:\on\to\on$ be continuous, 
and $U=Q_\al$, $V=Q_\ba$ 
be trees in $\jqo$. 
To prove that $\huf UfV$ weakly covers $U$,    
let $h=f_{\al\ba}$. 
Then $\huf QhQ$ weakly covers $Q$ by Step~4. 
Thus there is an ordinal $\xi<\omi$ 
such that if $q\in Q$ and $q(j)\ge\xi$, 
$j<\lh q$ then there is $n_0=n_0(q)\ge\lh q$ satisfying 
$s\in \huf QhQ$ whenever $s\in\brn Qn$, $n\ge n_0$, 
$q\sq s$. 
We claim that $\xi$ witnesses that $\huf UfV$ 
weakly covers $U$. 

Assume that $u\in U=Q_\al$, $u(j)\ge\xi$, $j<\lh u$. 
Then $q=\al\we u\in Q$. 
Let $n_0=n_0(q)$ and $m_0=n_0-1$. 
Now suppose that $m\ge m_0$ and $u\sq t\in\brn Um$; 
\rit{it remains to prove that\/ $t\in \huf UfV$}. 
Let $n=m+1$ and $q=\al\we u$, $s=\al\we t$. 
Thus $n\ge n_0$ and $q\sq s\in\brn Qn$, 
so that $s\in \huf QhQ$, by the choice of $n_0$. 

{\it Case 1\/}: 
$[Q]\cap (\ima h{[\req Qs]})$ is bounded. 
However $h=f_{\al\ba}$ and $U=Q_\al$, $V=Q_\ba$, hence 
$[Q]\cap (\ima h{[\req Qs]})=
\ba\we ([V]\cap(\ima f{[\req Ut]}))$. 
Thus the set $[V]\cap(\ima f{[\req Ut]})$ is bounded, 
therefore $t\in \huf UfV$.

{\it Case 2\/}: 
$h\res {[\req Qs]}$ is a total identity,   
$h(x)=x$ whenever $x\in \req Qs$. 
Then $\ba=\al$, $U=V$, and 
$f\res {[\req Ut]}$ is a total identity,    
thus still $t\in \huf UfV$.
}

\parf{Blocking sequences and the forcing} 
\las{vop}

{\ubf We argue in the constructible universe $\rL$} 
in this section.

The forcing to prove Theorem~\ref{Tun} 
will be defined as the union of a \dd\omi sequence  
of \sfo s of size $\ali$, increasing in the sense of 
a relation $\prok$  
(Definition~\ref{dfex}). 
We here introduce the notational system to be used 
in this construction.

\bdf
\lam{vstf}
Let $\stf$ be the set of all \sfo s of 
cardinality $\le\ali$. 
\imar{stf}%
\index{forcing!$\stf$}%
\index{zzWTF@$\stf$}%

If $\vjpi=\sis{\jpi_\al}{\al<\la}$ 
is a sequence of forcings $\jpi_\al\in\stf$, 
then let $\bigcup\vjpi=\bigcup_{\al<\la}\jpi_\al$, 
\index{zzbigcupP@$\bigcup\vjpi$}%
and let $\cM(\vjpi)$ be the least transitive model  
\index{zzMP@$\cM(\vjpi)$}%
\index{model!zzMP@$\cM(\vjpi)$}%
of $\zfcm$ of the form $\rL_\vt$, containing $\vjpi$, 
in which both $\la$ and the set 
$\bigcup\vjpi$ are of cardinality $\le\ali$.

If $\la\le\omd$ then let $\vstf_\la$ 
be the set of all \dd\la sequences    
$\vjpi=\sis{\jpi_\al}{\al<\la}$ 
of forcings $\jpi_\al\in\stf$, 
satisfying the following:
\ben
\fenu
\itlb{vstf*} 
if\, $\ga<\la$\, then\, 
${\bigcup{(\vjpi\res\ga)}}\prol{\cM(\vjpi\res\ga)} 
\jpi_\ga$. 
\een
Let $\vstf=\bigcup_{\la<\omd}\vstf_\la$.
\edf

The set $\vstf\cup\vsto$ is ordered by the 
extension relations $\su$ and $\sq$. 

\ble
\lam{202}
Assume that\/ $\ka<\la<\omd$, 
and\/ $\vjpi=\sis{\jpi_\al}{\al<\ka}\in\vstf$. 
Then$:$
\ben
\renu
\itlb{2021}%
the union\/ $\jpi=\bigcup\vjpi$ belongs to\/ 
$\stf\,;$

\itlb{2022}%
there is a sequence\/ $\vjqo\in\vstf$ 
such that\/ $\dom(\vjqo)=\la$ and\/   
$\vjpi\su\vjqo\;.$
\een
\ele
\bpf
To prove \ref{2022} apply Theorem~\ref{tex} by induction 
on $\la$.
\epf

\bdf
[key definition]
\lam{dzap}
A sequence\/ $\vjpi\in\vstf$
\rit{blocks} a set $W\sq\vstf$ if either  
\index{blocks}%
$\vjpi\in W$ 
or there is no sequence 
$\vjqo\in W$ satisfying $\vjpi\sq\vjqo$. 
\edf 

{\ubf Sets $\hka\ka$ and definability classes.}
\imar{hka ka}
Recall that $\hka\ka$ is the set of all 
\index{hereditarily of cardinality $<\ka$!$\hka\ka$}%
\index{zzhk@$\hka\ka$}%
sets \rit{hereditarily of cardinality $<\ka$}.
Thus $x\in\hka\ka$ if the \rit{transitive closure} $\TC(x)$ 
is a set of cardinality $<\ka$.
\imar{hcd}
In particular 
$\hc=\hka{\omi}$ is the set of all 
\rit{hereditarily (at most) countable}
\index{hereditarily countable!$\hc$}%
\index{zzhc@$\hc$}%
sets,  while 
$\hcd$ is the set of all 
\index{hereditarily of cardinality $\le\ald$!$\hcd$}%
\index{zzhc2@$\hcd$}%
sets \rit{hereditarily of cardinality ${\le\ali}$};  
$\hc=\rL_{\omi}$ and $\hcd=\rL_{\omd}$ in 
the constructible universe $\rL$.

\imar{hs ka n}%
\imar{gs ka n}%
$\hs\ka n$, resp., $\gs\ka n$ is the class 
of all sets $X\sq\hka\ka$, definable in $\hka\ka$ by 
a $\is{}n$ formula with parameters in $\hka\ka$, 
resp., with no parameters. 
The classes $\hp\ka n$, $\gp\ka n$ have the same meaning 
(with $\ip{}n$ formulas), 
and $\hd\ka n=\hs\ka n\cap\hp\ka n$, 
$\gd\ka n=\gs\ka n\cap\gp\ka n$, 
as usual. 
In particular, $\hd\ka 0=\hs\ka 0=\hp\ka 0$ and 
$\gd\ka 0=\gs\ka 0=\gp\ka 0$ 
(definability by bounded formulas, with/without parameters).
See more on \dd\in definability 
in \cite[Part B, Chap.\,5,\;Sect.\,4]{skml} or elsewhere.%
\index{definability classes!k1@$\gs\ka n,\,\gp\ka n,\,\gd\ka n$}%
\index{zzSHCk1@$\gs\ka n,\,\gp\ka n,\,\gd\ka n$}%
\index{definability classes!k2@$\hs\ka n,\,\hp\ka n,\,\hd\ka n$}%
\index{zzSHCk2@$\hs\ka n,\,\hp\ka n,\,\hd\ka n$}%
%

In particular, we consider the classes 
\imar{ags n}
$\ags n,\,\agp n,\,\agd n$ of definability 
in $\hcd$ (parameters not allowed) and 
\imar{ahs n}
$\ahs n,\,\ahp n,\,\ahd n$ 
(all parameters in $\hcd$ allowed) --- 
this is the case $\ka=\ald$ in the above definitions.
\index{definability classes!$\ags n,\,\agp n,\,\agd n$}%
\index{zzSHChc@$\ags n,\,\agp n,\,\agd n$}%
\index{definability classes!$\ahs n,\,\ahp n,\,\ahd n$}%
\index{zzSHChc@$\ahs n,\,\ahp n,\,\ahd n$}%

\bte
[the blocking sequence theorem, in $\rL$]
\lam{ep}
Let\/ $\nn\ge2$. 
There exists a sequence\/ 
$\vdp=\sis{\dP_\al}{\al<\omd}\in\vsto$  
satisfying the following two conditions$:$
\ben
\renu
\itlb{ep1}%
$\vdp$, 
as the set of pairs\/ $\ang{\al,\dP_\al}$, 
belongs to the definability class\/ $\agd{\nn-1}\;;$

\itlb{ep2}%
if\/ $\nn\ge3$ and\/ $W\sq\vstf$ is a\/ $\ahs{\nn-2}$ set\/ 
%
then there is an ordinal\/ $\ga<\omd$ such that 
the restricted sequence\/ 
$\vdp\res\ga=\sis{\dP_\al}{\al<\ga}\in\vstf$
blocks\/ $W$.
\een
\ete
\bpf
Let $\lel$ be the canonical $\id{}1$ wellordering of $\rL$; 
thus its restriction to $\hcd=\rL_{\omd}$ is $\agd1$.
As $\nn\ge3$, there exists a universal $\ags{\nn-2}$ set 
$\gU^\nn\sq\omd\ti\hcd$. 
That is, $\gU^\nn$ is $\ags{\nn-2}$ 
(parameter-free $\is{}{\nn-2}$ definable in $\hcd$), 
and for every set $X\sq\hcd$ of type  
$\ahs{\nn-2}$ 
($\is{}{\nn-2}$ definable in $\hcd$ with arbitrary parameters) 
there is an ordinal $\da<\omi$ such that 
$X=\gU^\nn_\da=\ens{x}{\ang{\da,x}\in\gU^\nn}$. 
The choice of $\omd$ as the domain of parameters is 
validated by the assumption $\rV=\rL$, which implies the 
existence of a $\agd1$ surjection $\omd\onto\hcd$.

Coming back to Definition~\ref{dzap}, note that for any 
sequence $\vjpi\in\vstf$ and any set $W\sq\vstf$ 
there is a sequence $\vjqo\in\vstf$ which satisfies 
$\vjpi\su\vjqo$ and blocks $W$.
This allows us to define $\vjqo_\al\in\vstf$ by induction 
on $\al<\omi$ so that $\vjqo_0=\pu$, 
$\vjqo_\la=\bigcup_{\al<\la}\vjqo_\al$, and each 
$\vjqo_{\al+1}$ is equal to the \dd\lel least 
sequence $\vjqo\in\vstf$ which satisfies 
1) $\vjqo_\al\su\vjqo$ and 2)  
if $\nn\ge3$ then $\vjqo$ blocks $\gU^\nn_\al$. 

Then $\vdp=\bigcup_{\al<\omd}\vjqo_\al\in\vsto$. 
Condition \ref{ep2} holds by construction, while  
\ref{ep1} follows by a routine verification, 
based on the fact that $\vstf\in\agd1$ and 
$\gU^\nn\in \ags{\nn-2}$ (provided $\nn\ge3$).
\epf

\bdf
[in $\rL$]
\lam{vdp}
We fix a sequence  
$\vdp=\sis{\dP_\al}{\al<\omd}\in\vsto$, 
given by Theorem~\ref{ep} for a number 
$\nn\ge2$, for which  
\index{zzn@$\nn\ge 2$}%
\index{n@number $\nn\ge 2$}%
Theorem~\ref{Tun} is to be established.

In particular $\vdp$ satisfies \ref{ep1} and \ref{ep2}
of Theorem~\ref{ep}.

If $\ga<\omd$ then let 
$\cM_\ga=\cM(\vdp\res\ga)$,   
\imar{vdp, mdp ga}
\index{zzMga@$\cM_\ga$}%
\index{model!zzMga@$\cM_\ga$}%
and 
$\mdp{\ga}=\bigcup_{\al<\ga}\dP_\al$, 
$\dP=\bigcup_{\al<\omd}\dP_\al$. 
\index{zzPd@$\dP$}%
\index{forcing!zzPd@$\dP$}%
\edf

\parf{Some forcing properties} 
\las{pres+}

The \sfo\ $\dP\in\rL$ defined by \ref{vdp} 
will be the forcing notion for the proof 
of Theorem~\ref{Tun}. 
The next lemma establishes some properties of  
$\dP$.

{\ubf We continue to argue in $\rL$ in the conditions and 
notation of Definition~\ref{vdp}}.

\ble
\lam{jden}
$\dP$ is a\/ \sfo, 
all sets\/ $\dP_\al$, $\mdp\ga$ belong to\/ $\stf$. 
In addition$:$
\ben
\renu
\itlb{prop1}%
if\/ $\al<\omd$ then\/ ${\mdp\ga}\prol{\cM_\ga}\dP_\ga\,;$

\itlb{jden1}%
if\/ $\al<\omd$ and the set\/ 
$D\in\cM_\al\yt D\sq\mdp\al$ 
is pre-dense in\/ $\mdp\al$ then it is 
pre-dense in\/ $\dP$, too$;$ 

\itlb{jden2}%
every set\/ $\dP_\al$ is pre-dense in\/ $\dP;$

\itlb{jden33}%
if\/ $\al<\omd$ and trees\/ $T,T'\in\mdp{\al}$   
are incompatible in\/ $\mdp{\al}$ then\/ $T,T'$  
are incompatible in\/ $\dP$, too$;$ 

\itlb{prop2}%
if\/ $f:\on\to\on$ is continuous then the set of 
all trees\/ $T\in\dP$ such that\/ $\ima f{[T]}$ is 
bounded or\/ $f\res[T]$ is a bijection, is dense 
in\/ $\dP\,;$  

\itlb{prop4}%
if\/ $f:\on\to\on$ is continuous then the set of 
all trees\/ $T\in\dP$ such that
$(1)$ 
$f\res[T]$ is a total identity, {\bf or}, 
for some\/ $\ga<\omd$, 
$(2)$ 
$f\res {[T]}$\/ {\ubf avoids} 
\index{avoids}
$\dP_\ga$ in the sense that if\/ $V\in\dP_\ga$ 
then the subset\/ 
$\ens{s\in T}{[V]\cap(\ima f{[\req Ts]}\text{ is bounded}}$
weakly covers\/ $T$, 
is dense in\/ $\dP\,;$  

\itlb{prop3}%
if\/ $\nn\ge3$ 
and a set\/ $Q\sq\pew$ belongs to\/ $\ahs{\nn-2}$, 
then\/ $\dP\cap{(Q\cup Q^-)}$ is dense in\/ $\dP$, 
where\/ $Q^-=\ens{T\in\pew}{\neg\:\sus S\in Q\,(S\sq T)}$.
\een
\ele

\bpf
\ref{prop1}
holds by \ref{vstf*} of Definition~\ref{vstf}.\pagebreak[1]

\ref{jden1}  
We use induction on $\ga\yd \al\le\ga<\omd$, 
to check that if $D$ is pre-dense in $\mdp\ga$  
then it remains pre-dense in   
$\mdp\ga\cup\dP_\ga=\mdp{\ga+1}$ by \ref{prop1}
and Lemma~\ref{tec}\ref{tec1}. 
Limit steps, including the final step to $\dP$ 
($\ga=\omd$) are routine. 

\ref{jden2}
$\dP_\al$ is dense in  
$\mdp{\al+1}=\mdp\al\cup\dP_\al$ 
by \ref{dfex}\ref{dfex1}.
It remains to refer to \ref{jden1}.

\ref{jden33}
Prove by induction on $\ga$ that if $\al<\ga\le\omi$ 
then $T,T'$ are incompatible in $\mdp{\ga}$, using 
\ref{prop1} and Lemma~\ref{tec}\ref{tec2}. 

To prove \ref{prop2} and \ref{prop4}
let $T_0\in\dP$.
There is an ordinal $\ga<\omd$ such that 
$T_0\in\mdp\ga$ and $f\in\cM_\ga$. 
We have ${\mdp\ga}\prol{\cM_\ga}\dP_\ga$ by \ref{prop1}. 
Therefore by \ref{dfex5} of Definition~\ref{dfex} there 
is a tree $T\in \dP_\ga$ such that $T\sq T_0$ and 
$\ima f{[T]}$ is bounded or\/ $f\res[T]$ is a bijection, 
so we get \ref{prop2}.
Further by \ref{dfex6} of Definition~\ref{dfex} 
if $V\in\dP_\ga$ then 
the set $\huf TfV$, of all strings $s\in T$  
such that $[V]\cap(\ima f{[\req Ts]})$ is bounded or 
$f\res[\req Ts]$ is a total identity, weakly covers $T$. 
We have two cases.

{\ubf Case 1\/}: 
$f\res[\req Ts]$ is a total identity for at least one 
$s\in T$. 
Then the corresponding subtree $T'=\req Ts$ 
satisfies (1) of \ref{prop4}.

{\ubf Case 2\/}: 
for each $V\in\dP_\ga$,   
the set $H(V)$ of all strings $s\in T$  
such that $[V]\cap(\ima f{[\req Ts]})$ is bounded, 
weakly covers $T$, thus $T$ itself 
satisfies (2) of \ref{prop4}. \vom

\ref{prop3}
Suppose that $\nn\ge3$. 
Let $T_0\in\dP$, that is, 
$T_0\in\mdp{\al_0}$, $\al_0<\omd$. 
The set $W$ of all sequences $\vjpi\in\vstf$, 
such that $\vdp\res\al_0\sq\vjpi$ and 
$\sus T\in Q\cap(\bigcup\vjpi)\,(T\sq T_0)$, 
belongs to $\ahs{\nn-2}$ 
along with $Q$. 
Therefore there is an ordinal $\al<\omd$ such that 
$\vdp\res\al$ blocks $W$. 
We have two cases.

{\ubf Case 1:\/} 
$\vdp\res\al\in W$.
Then 
the related tree $T\sq T_0$ belongs to 
$Q\cap\dP$.

{\ubf Case 2:\/} 
there is no sequence in $W$ which extends  
$\vdp\res\al$.
Let $\ga=\tmax\ans{\al,\al_0}$. 
Then $\mdp\ga\prol{\cM_\ga}\dP_\ga$ by \ref{prop1}. 
As $\al_0\le\ga$, there is a tree $T\in\dP_\ga$, 
$T\sq T_0$. 
We claim that $T\in Q^-$, which completes the 
proof in Case 2. 

Suppose to the contrary that $T\nin Q^-$, 
thus there is a tree $S\in Q$, $S\sq T$.
The set 
$\bR=\dP_\ga\cup\ens{\req St}{t\in S}$ 
is a \sfo\ and obviously
$\dP_\ga\prok\bR$, hence 
still $\mdp\ga\prol{\cM_\ga}\bR$ holds by 
Lemma~\ref{teb}. 
It follows that the sequence $\vjR$ defined by  
$\dom\vjR=\ga+1$, 
$\vjR\res\ga=\vdp\res\ga$, and $\vjR(\ga)=\bR$, 
belongs to $\vstf$, and even $\vjR\in W$ since 
$S\in Q\cap\bR$.  
Yet $\vdp\res\al\su\vjR$,
which contradicts to the Case 2 hypothesis. 
\epf


To prove a chain condition for $\dP$, 
we'll need the following general lemma.
See Definition~\ref{vdp} on models $\cM_\al$.

\ble
[in $\rL$]
\lam{club}
If\/ $X\sq\hcd=\rL_{\omd}$ then the set\/ $\skri O_X$ 
of all ordinals\/ $\al<\omd$ such that the model\/ 
$\stk{\rL_\al}{X\cap\rL_\al}$  
is an elementary submodel of\/ $\stk{\rL_{\omd}}{X}$ 
and\/ $X\cap\rL_\al\in\cM_\al$, is unbounded in\/ $\omd$.
\vyk{
Generally if\/ $X_n\sq\hcd$ for all\/ $n$  
then the set\/ $\skri O$ of all ordinals\/ 
$\al<\omd$ such that\/  
$\stk{\rL_\al}{\sis{X_n\cap\rL_\al}{n<\om}}$ is  
an elementary submodel of\/  
$\stk{\rL_{\omi}}{\sis{X_n}{n<\om}}$ 
and\/ $\sis{X_n\cap\rL_\al}{n<\om}\in\cM_\al$, 
is unbounded in\/ $\omd$.
}%
\ele
\bpf
Let $\al_0<\omd$. 
There is an elementary submodel $M$ 
of $\stk{\rL_{\om_3}}{\in}$, of cardinality 
$\card M=\ali$,  
which contains $\al_0\yi\omd\yi X$  
and is such that the set $M\cap\rL_{\omd}$ is transitive. 
Consider the Mostowski collapse $\phi:M\onto\rL_\la$.  
Let $\al=\phi(\omd)$. 
Then $\al_0<\al<\la<\omd$ and $\phi(X)=X\cap\rL_\al$ 
by the choice of $M$. 
We conclude that $\stk{\rL_\al}{X\cap\rL_\al}$ is  
an elementary submodel of $\stk{\rL_{\omd}}{X}$.
And $\card\al>\ali$ in $\rL_\la$, hence  
$\rL_\la\sq\cM_\al$. 
Then $X\cap\rL_\al\in\cM_\al$, as   
$X\cap\rL_\al\in\rL_\la$ by construction.
%
\epf

\bcor
[in $\rL$]
\lam{ccc}
{\rm(i)}
If\/ $A\sq\dP$ is an antichain then\/ $\card A\le\ali$.
\ben
\renu
\atc
\itlb{ccc1}%
Let\/ $D_n\sq\dP$ be pre-dense in\/ $\dP$, 
for each\/ $n$. 
Then the set of all trees\/ $T\in\dP$, 
satisfying\/ 
$\kaz n\,(\brn Tn\sq \hdt{D_n}T)$, 
is dense in\/ $\dP$.
\een
\ecor
\bpf
(i)
Let $A\sq\dP$ be a maximal antichain. 
By Lemma~\ref{club} there is an ordinal $\al$ 
such that $\stk{\rL_\al}{\dP',A'}$ is  
an elementary submodel of $\stk{\rL_{\omd}}{\dP,A}$, 
where $\dP'=\dP\cap\rL_\al$ and $A'=A\cap\mdp\al$, 
and in addition $\dP',A'\in\cM_\al$. 
By the elementarity, we have  
$\dP'=\mdp\al$ and $A'=A\cap\mdp\al\in\cM_\al$, 
and $A'$ is a maximal antichain, hence a pre-dense  
set, in $\mdp\al$. 
But then $A'$ is a pre-dense set, hence, 
a maximal antichain, in the whole set $\dP$ 
by Lemma~\ref{jden}\ref{jden1}. 
Thus $A=A'$, and $\card A=\card{A'}\le\ali$.

\ref{ccc1}
We wlog assume that all $D_n$ are open dense, for if not 
then replace $D_n$ by the set
$\ens{T\in\dP}{\sus S\in D_n\,(T\sq S)}$. 
Let $T_0\in\dP$. 
Pick a maximal antichain $A_n\sq D_n$ in each $D_n$. 
Then all sets $A_n$ are maximal antichains in $\dP$ 
by the open density, and $\card{A_n}\le\ali$ by 
(i). 
Therefore there is an ordinal $\al<\omd$ such that 
the set $A=\bigcup_nA_n$ satisfies $A\sq\mdp\al$
and $A$, $T_0$, and 
the sequence $\sis{A_{n}}{n<\om}$ belong to $\cM_\al$. 
By the maximality of $D_n$ and Lemma~\ref{jden}\ref{jden33},
each $D'_n=D_n\cap\mdp\al$ is 
dense in $\mdp\al$. 
It follows by Lemma~\ref{jden}\ref{prop1} and 
\ref{dfex4} of Definition~\ref{dfex} that there is a tree 
$T\in\dP_\al$ such that $T\sq T_0$ and 
$\brn Tn\sq \hdt{D_n}T$
for all $n$.
\epf

\vyk{
\bcor
[in $\rL$]
\lam{Nden}
Let\/ $D_n\sq\dP$ be pre-dense in\/ $\dP$, 
for each\/ $n$. 
Then the set of all trees\/ $T\in\dP$, 
satisfying\/ 
$\kaz n\,\kaz s\in\brn Tn\,(\req Ts\ine D_n)$, 
is dense in\/ $\dP$.
\ecor
\bpf

\epf
}

\parf{The model} 
\las{bex1}

This section presents some key properties of  
\dd\dP generic extensions $\rL[G]$ of $\rL$ 
obtained by adjoining a \dd\dP generic  
set $G\sq\dP$ to $\rL$. 
Recall that the forcing notion $\dP\in\rL$ 
was introduced by Definition~\ref{vdp}, along with 
some related notation. 

\bcor
\lam{calg}
If a set\/  
$G\sq\dP$ is\/ \dd\dP generic over\/ $\rL$  
then\/ $\omil<\omi^{\rL[G]}=\omdl$. 
\ecor
\bpf
That $\omil<\omi^{\rL[G]}$ follows from the fact
that $\xx G$ is a cofinal map $\om\to\omil$.
To prove $\omi^{\rL[G]}=\omdl$
use Corollary~\ref{ccc}.
\epf

\bgg
\lam{ba1}
Arguing in generic extensions of $\rL$, we'll use
standard notation like $\om_\xi^\rL$ to denote
\dd\rL cardinals.
We also use $\pel$ to denote
``the set $\pew$ defined in $\rL$''.
Thus for instance $\dP\sq\pel$.
\egg

We'll make use of a 
{\ubf coding system for continuous maps}, 
helpful whenever  
``the same'' continuous $f:\on\to\oom$ 
is considered in different models.

\bdf
\lam{dccf}
Let $\vt\in\Ord$.
A \rit{code of continuous function}
\index{code of continuous function, $\ccf_\vt$}%
\index{zzccfthet@$\ccf_\vt$}%
from $\ol$ to $\vt^\om$ is any 
map $\rc:\dom\rc\to\vt$ with $\dom\rc\sq\bsdl\ti\om$, 
such that the sets 
\index{zzScnxi@$\yc\rc n\xi$}%
$\yc\rc n\xi=
\ens{s\in\bsdl}{\ang{s,n}\in\dom\rc\land\rc(s,n)=\xi}$
satisfy the following for any $n$:\vtm

(1) 
if $\xi\ne\et$, $u\in \yc\rc n\xi$, $v\in \yc\rc n\et$, then 
$u,\,v$ are
incompatible, and\vom 

(2) 
\index{zzKcn@$\yc\rc n{}$}%
$\yc\rc n{}=\bigcup_\xi\yc\rc n\xi$ 
is a belt for $\bsdl$, 
\ie, 
$\kaz x\in\ol\,\sus m\,(x\res m\in\yc\rc n{})$.\vtm 

\noi
Let $\ccf_\vt$ be the set of all such codes.
If $\rc\in\ccf_\vt$ then a continuous
$f_\rc:\ol\to\vt^\om$ 
\index{zzfc@$f_\rc$}
is defined as follows. 
If $x\in\ol$ and $n<\om$, then by definition
there is a unique $\xi<\vt$ such that
$x\res k\in \yc\rc n\xi$ for some 
$k$.
Let $f_\rc(x)(n)=\xi$.
\edf

If $f:\ol\to\vt^\om$ is continuous then its code
\index{zzcodef@$\koe f$}%
\index{code!$\koe f$}%
$\rc=\koe f\in\ccf_\vt$ is defined by
$\yc\rc n\xi=
\ens{s\in\bsdl}{\kaz x\in\ol\,(s\su x\imp f(x)(n)=\xi)}$;
then $f_\rc=f$.

\bre
[absoluteness]
\lam{rccf}
Being a code in $\ccf_\vt$ is absolute since so is the condition 
of being a belt, see Remark~\ref{belta}.                  
\ere

\vyk{
If $G\sq\dP$ is generic over $\rL$ and 
$\rc\in\ccf_\vt$ 
then $\xx G\in\ol$, hence we can define 
$\rc[G]=f_\rc(\xx G)\in\vt^\om$; 
the definition of $\xx G$ see Remark~\ref{Pfo}.
}

\ble
\lam{repd}
If\/ $G\sq\dP$ is generic over\/ $\rL$,
$\vt\in\Ord$, $y\in\vt^\om\cap\rL[G]$, 
then
\ben
\renu
\itlb{repd1}%
there is a code\/ $\rc\in\ccf_\vt\cap\rL$  
such that\/ $y=f_\rc(\xx G)\,;$ 

\itlb{repd2}%
if\/ $\vt=\omil$ then\/
$y$ is bounded in\/ $\omil$ or\/ $G\in\rL[y]\,;$

\itlb{repd3}%
if\/ $\vt=\omil$ and\/ 
$y$ is unbounded in\/ $\omil$ then\/ $y=\xx G$ 
or there is an ordinal\/ $\ga<\omd^\rL$ such that\/ 
$y\nin\bigcup_{V\in\dP_\ga}[V]$. 
\een
\ele
\bpf
\ref{repd1} 
There is a \dd\dP name $t\in\rL$   
satisfying $y=t[G]$ 
(the \dd Gvaluation of $t$). 
It can be assumed that $\dP$ forces that $t$   
is valuated as an element of $\vt^\om$. 

{\ubf Arguing in $\rL$}, 
let  
$\tau_{n\xi}=\ens{T\in\dP}{T\text{ forces }t(n)=\xi}$ 
($n<\om$ and $\xi<\vt$). 
The sets $\tau_n= \bigcup_\xi \tau_{n\xi}$ 
are open dense in $\dP$. 
It follows by Corollary~\ref{ccc}\ref{ccc1} that there is a 
tree $T\in G$ such that $\req Ts\in \tau_n$  
whenever $n<\om$ and $s\in\brn Tn$.
This allows us to define,
{\ubf still in $\rL$}, a continuous
$f':[T]\to\vt^\om$ by $f'(x)(n)=\xi$ iff
the only string $s\in\brn Tn$ with $s\su x$
belongs to $\tau_{n\xi}$.
Let $f:\on\to\vt^\om$ be a continuous extension
of $f'$.
Then $\rc=\koe f\in\ccf_{\vt}\cap\rL$, 
and easily $y=f_\rc(\xx G)$.

\ref{repd2} 
Let, by \ref{repd1}, $\rc\in\ccf_{\omil}\cap\rL$    
and $y=f_\rc(\xx G)$.
By Lemma~\ref{jden}\ref{prop2}, there is a tree 
$T\in G$ such that, in $\rL$, $\ima {f_\rc}{[T]}$ 
is bounded or $f_\rc\res [T]$ is a bijection. 

{\it Case 1\/}: 
in $\rL$, $\ima {f_\rc}{[T]}$ is bounded, that is, 
there is an ordinal $\ba<\omil$ satisfying 
$f_\rc({x})\in\ba^\om$ for all ${x}\in[T]\cap\rL$. 
But $f_\rc$ is continuous while 
$[T]\cap\rL$ is dense in $[T]$ in $\rL[G]$. 
It follows that $f_\rc({x})\in\ba^\om$ for all 
${x}\in[T]\cap \rL[G]$. 
In particular $y=f_\rc(\xx G)\in\ba^\om$ since 
$\xx G\in [T]$  (because $T\in G$), so $y$ is 
bounded.

{\it Case 2\/}: 
in $\rL$,  $f_\rc\res [T]$ is a bijection. 
The bijectivity is equivalent to the wellfoundedness 
of the tree $W_\rc$ of all pairs $\ang{s,t}$ of strings 
$s,t\in T$ such that $\lh s=\lh t$ and there exist 
no strings $u,v$ satisfying: $u\sq s$, $v\sq t$, 
and $u\in\yc\rc n\xi$, $v\in\yc\rc n\et$ for some $n$ 
and $\xi\ne\et$.
Therefore the bijectivity of $f_\rc\res [T]$ is an 
absolute property of $\rc,T$. 
Thus $f_\rc\res [T]$ is a bijection in $\rL[G]$, and 
we have $\xx G= f_\rc\obr(y)\in\rL[y]$, as required.

\ref{repd3}
We still assume that, by \ref{repd1}, $y=f_\rc(\xx G)$, 
where $\rc\in\ccf_{\omil}\cap\rL$.
By Lemma~\ref{jden}\ref{prop4}, there is a tree 
$T\in G$ such that, in $\rL$, ${f_\rc}\res{[T]}$ 
is a total identity or, for some $\ga<\omdl$, 
$f_\rc\res [T]$ avoids $\dP_\ga$ in the sense of 
\ref{jden}\ref{prop4}. 

{\it Case 1\/}: 
in $\rL$, ${f_\rc}\res{[T]}$ is a total identity, 
that is, $f_\rc({x})={x}$ for all ${x}\in[T]\cap\rL$. 
By the same simple continuity/density argument, 
we have $f_\rc({x})={x}$ for all 
${x}\in[T]\cap \rL[G]$, 
in particular $y=f_\rc(\xx G)=\xx G$.

{\it Case 2\/}: 
$\ga<\omdl$ and, in $\rL$: $f_\rc\res [T]$ avoids $\dP_\ga$, 
that is, if\/ $V\in\dP_\ga$ 
then the subset\/ 
$T(V)=\ens{s\in T}
{[V]\cap(\ima f{[\req Ts]})\text{ is bounded}}$
(defined in $\rL$) weakly covers\/ $T$. 
Now let $V\in\dP_\ga$ and check that $y\nin [V]$. 
By the Case 2 assumption, $T(V)$ weakly covers $T$. 
Therefore, as $\xx G\in[T]$ is definitely unbounded, 
there is a string $s\in T(V)$ satisfying $s\sq\xx G$. 
Then $S=\req Ts\in G$ and $[V]\cap(\ima f{[S]})$ is bounded, 
so that there is an ordinal $\ba<\omil$ satisfying: 
if ${x}\in[S]\cap\rL$ and $f_\rc({x})\in[V]$ 
then $f_\rc({x})\in\ba^\om$. 
We claim that 
\rit{the implication\/ 
${f_\rc({x})\in[V]}\imp {f_\rc({x})\in\ba^\om}$
also holds for all\/ ${x}\in[S]\cap\rL[G]$.} 
Assume that this is established. 
As ${x}=\xx G\in[S]$ (because $S\in G$), we then have 
${y\in[V]}\imp {y\in\ba^\om}$. 
(Recall that $y=f_\rc(\xx G)$.) 
But $y$ is unbounded, hence $y\nin[V]$, 
as required.

To prove the claim, let 
${x}_0\in[S]\cap\rL[G]$ be a counterexample, so 
$y_0=f_\rc({x}_0)\in[V]$ but 
$y_0(n_0)=\xi$ for some $n_0$ and $\xi\ge\ba$.
The existence of such ${x}_0$ is equivalent to the 
non-wellfoundedness 
of the tree $W$ of all strings $s\in S$ such that 
$s\in\yc\rc{n_0}\et$ for all $\et\ne\xi$,  
and there is no string $u\nin V$ satisfying: 
$\kaz j<\lh u\,(s\in\yc\rc j{u(j)})$.
Therefore the existence of ${x}_0$ is an 
absolute property of $\rc,S,V$. 
Thus such an ${x}_0\in[S]$ exists already in $\rL$, 
contrary to the Case 2 assumption.
\epf

\bcor
\lam{dfs}
Let\/ $G\sq\dP$ be generic over\/ $\rL$. 
Then it is true in\/ $\rL[G]$ that\/
\ben
\renu
\itlb{dfs1}%
$\xx G$ 
is the only member of the intersection\/ 
$\bigcap_{\ga<\omi}\bigcup_{T\in\dP_\ga}[T]\;;$  

\itlb{dfs2}%
$\ans{\xx G}$ is a\/ $\ip\hc{\nn-1}$ singleton$;$

\itlb{dfs3}%
there is a\/ $\ip1{\nn}$ real singleton\/ $\ans{r}$,
$r\in\bn$, such that\/ $\rL[r]=\rL[\xx G]$.
\een
\ecor
\bpf
\ref{dfs1}
Each $\dP_\ga$ is pre-dense in $\dP$ by
Lemma~\ref{jden}\ref{jden2}.
It follows that $\xx G\in \bigcup_{T\in\dP_\ga}[T]$,
by the genericity.
The uniqueness follows from Lemma~\ref{repd}\ref{repd3}.

\ref{dfs2}
The sequence
$\vdp= 
\ens{\ang{\ga,\dP_\ga}}{\ga<\omdl}$
is of type $\agd{\nn-1}$ in $\rL$ by
Definition~\ref{vdp}.
However $\hcd$ in the sense of $\rL$ coincides with
the constructible part of $\HC$
(= hereditarily countable sets)
in the sense of $\rL[G]$, because $\omi^{\rL[G]}=\omdl$
by Corollary~\ref{calg}.
It easily follows that $\vdp$
is $\id\hc{\nn-1}$ in $\rL[G]$.
On the other hand,
$$
\ans{\xx G}=\ens{x}{\kaz\ga\,\kaz\dQ\,
(\ang{\ga,\dQ}\in\vdp \imp \sus T\in\dQ\,(x\in[T])}
$$
by \ref{dfs1}. 
This yields the result since $\sus T\in\dQ$ is a
bounded quantifier.

\ref{dfs3}
If $r\in\bn$ then let $(r)_n(k)=r(2^n(2k+1)-1)$, thus
$(r)_n\in\bn$.
Let $W$ be the $\ip11$ set of all reals which code an
ordinal, and let $|w|<\omi$ be the ordinal coded
by $w\in W$.
Let $r\in\bn$ be defined so that each $(r)_n$ belongs
to $W\cap\rL$ and is \dd\lel minimal of all
$w\in W\cap\rL$ satisfying $|w|=\xx G(n)$.
Thus $r$ is a real in $\rL[G]$.
The singleton $\ans r$ is defined in $\hc$ of $\rL[G]$
by the following formula:
\begin{quote}
$\kaz n$, $(r)_n\in W\cap\rL$ and
$(r)_n$ is \dd\lel minimal of all
$w\in W\cap\rL$ with $|w|=|(r)_n|$, and
$
\kaz x\in\oom\,(\kaz n\,(x(n)=|(r)_n|)\imp x=\xx G).
$
\end{quote}
It easily follows by the result of \ref{dfs2}
that $\ans r$ is a $\ip\hc{\nn-1}$ singleton as well,
hence a $\ip1{\nn}$ singleton.
\epf

Corollary~\ref{calg} and
Corollary~\ref{dfs}\ref{dfs2},\ref{dfs3} account
for items (i), \ref{Tun2}, \ref{Tun3} of
Theorem~\ref{Tun}.
Item \ref{Tun4} of the theorem is based on different
ideas related to claim \ref{prop3} of Lemma~\ref{jden}.
From now on we work towards this goal.

\parf{Shoenfield's transformation of $\is12$ formulas} 
\las{abs}

The following useful transformation of $\is12$ formulas
involves an idea in the proof of the  Shoenfield
absoluteness theorem. 

\bgg
\lam{pqr}
From now on $p,q,r$ denote reals in $\bn$.
\egg

\bte
\lam{abst}
Let\/ $\vpi(p_1,\dots,p_n)$
be a\/ $\is12$ formula of the form
$$
\left.
\bay{l}
\vpi(p_1,\dots,p_n)
\,:=\, 
\sus q\,\kaz r\,\sus m\,
R(q\res m,r\res m,p_1\res m,\dots,p_n\res m)\,,
\\[0.5ex]
\text{\rm where } R\sq(\nse)^{n+2}, \ R\in\rL, \
\text{\rm and }q,r,p_i \text{\rm\ are variables over } \bn
\eay
\right\},
\eqno(\ast)
$$                
and\/ 
$\vt\ge\alil$ a cardinal in\/ $\rL$.
Then there is a relation\/
$Q=Q_\vt(R)\sq \bst\ti(\nse)^{n+1}$, 
$Q\in\rL$,
such that\/ $Q$ is\/ $\gd\vt0(R)$
as a subset of\/ 
$\hka\vt$ in\/ $\rL$\snos
{Meaning that the equality 
$Q=\ens{w\in\hka\vt}{\psi(w)}$ holds in $\rL$, where 
$\psi$ is a bounded formula with $R$ as the only parameter.}%
, 
and it holds in any generic
extension\/ $M$ of\/ $\rL$ 
with\/ $\vt\ge\omi^M$ that
if\/ $p_1,\dots,p_n\in\bn$ 
then
$$
\vpi(p_1,\dots,p_n)
\leqv
{\sus\chi\in\vto\,\sus q\in\bn\,\kaz m\,
Q(\chi\res m,q\res m,p_1\res m,\dots,p_n\res m)}
\;.
$$
\ete
\bpf
$\vpi(p_1,\dots,p_n)$ is equivalent to 
$
\sus q\,(W_{q,p_1,\dots,p_n}\text{ is wellfounded})\,,
$
where
$$
W_{q,p_1,\dots,p_n}
=\ens{u\in\nse}
{\kaz j\le\lh u\,\neg\:
R(q\res j,u\res j,p_1\res j,\dots,p_n\res j)}\,, 
$$
hence --- in any universe $M$ as in the theorem
--- to the formula:
$$
\sus q\in\bn\,\sus f:W_{q,p_1,\dots,p_n}\to\vt\,
(f\text{ is order-preserving})\,.
$$
By ``order-preserving'' we mean: if
$u,v\in W_{q,p_1,\dots,p_n}$ then
$u\leq_{\text{LS}}v \eqv f(u)\le f(v)$, where
$\leq_{\text{LS}}$ is the Lusin -- Sierpinski
(= Kleene -- Brouwer) order on strings. 

Fix a recursive bijection $k\mto s_k:\om\onto\nse$,
with the inverse 
bijection $\num:\nse\to\om$, so that $s=s_{\num(s)}$.
We assume that $\lh s\le\num(s)$, $\kaz s$.
Let
$$
W^m_{q,p_1,\dots,p_n}=
\ens{s\in W_{q,p_1,\dots,p_n}}{\num(s)<m}\,,
$$
a finite set.
Then $\vpi(p_1,\dots,p_n)$ is equivalent to the formula
$$
\sus q\in\bn\,\sus\chi\in\vto\,\kaz m\,
(
\chi\circ\num
\text{ is order-preserving on }W^m_{q,p_1,\dots,p_n})\,.
$$
($\chi\circ\num$ is the superposition.)  
The subformula in brackets depends on $\chi\res m$ and 
$q\res m,p_1\res m,\dots,p_n\res m$ only. 
In other words, we have a relation
$Q=Q_\vt(R)\sq \bst\ti(\nse)^{n+1}$, 
still $Q\in\rL$,  
such that $\vpi(p_1,\dots,p_n)$ is equivalent
to the formula  
$$
\sus\chi\in\vto\,\sus q\in\bn\,\kaz m\,
Q(\chi\res m,q\res m,p_1\res m,\dots,p_n\res m)\,.
\eqno(\dag)
$$
%
Namely $Q$ contains all tuples 
$\ang{\sg,v,u_1,\dots,u_n}$ 
of strings $\sg\in\vt\lom$ and 
$v,u_i\in\nse$ 
of same length $\lh\sg=\lh v = \lh{u_i}=$ some $m$, 
such that the superposition $\sg\circ\num$
(defined on the set $S_m=\ens{s_j}{j<m}$)  
is order-preserving on the set
$$
W^m_{v,u_1,\dots,u_n}
=\ens{u\in S_m}
{\kaz j\le \lh u\,\neg\:
R(v\res j,u\res j,u_1\res j,\dots,u_n\res j)}\,. 
$$
To see that $Q$ is a $\gd\vt0(R)$ subset of $\hka\vt$, 
note first of all that $\vt=\Ord\cap\hka\vt$, which eliminates 
$\vt$ and $\vt\lom$ from the list of parameters. 
In the rest, we skip a routine verification of all elements of the 
definition of $Q$ being expressible by bounded formulas.
\epf

\parf{Auxiliary forcing relation} 
\las{fs}


Here we introduce a key tool for the proof of 
claim \ref{Tun4} of Theorem~\ref{Tun}. 
This is a forcing-like relation $\fo$. 
It is not explicitly connected with the forcing notion 
$\dP$ 
(but rather connected with the full 
wide tree forcing $\pew$), 
however it will be compatible with $\dP$ for  
formulas of certain quantifier complexity 
(Theorem~\ref{f-t}). 
The crucial advantage of $\fo$ will be its invariance  
a certain group of transformations  
(Lemma~\ref{inv}), 
a property that cannot be expected for $\dP$.
This will be the key argument in the proof 
of Theorem~\ref{Tun} below in Section~\ref{finar}.

\bgg
\lam{tet}
From now on, we 
let $\tet=\omdl$, so $\tet=\omd$ in $\rL$ but $\tet=\omi$ 
in \dd\dP generic extensions of $\rL$. 
\egg

{\ubf We argue in $\rL$.}
%
We consider a language $\cL$ whose elementary formulas,
called $\lis12$ 
(in spite that they are looking more like $\is11$), 
are those of the form
$$
\left.
\bay{l}
\vpi(p_1,\dots,p_n)
\,:=\, 
\sus\chi\in\tso\,\sus q\in\bn\,\kaz m\,
Q(\chi\res m,q\res m,p_1\res m,\dots,p_n\res m)\,,
\\[0.5ex]
\text{where $Q\in\rL$,
$Q\sq \tsd\ti(\nse)^{n+1}$, 
$Q$ is a $\ahd 0$ set,} \\[0.5ex]
\text{
and $q,p_i$ are variables over $\bn$.}
\eay
\right\}
\eqno(1)
$$
The dual class $\lip12$ consists of formulas
$$
\left.
\bay{l}
\vpi(p_1,\dots,p_n)
\,:=\, 
\kaz\chi\in\tso\,\kaz q\in\bn\,\sus m\,
Q(\chi\res m,q\res m,p_1\res m,\dots,p_n\res m)\,,
\\[0.5ex]
\text{with the same specifications.}
\eay
\right\}
\eqno(2)
$$
Higher classes $\lis1k$ and $\lip1k$ are
defined naturally, \eg\ $\lis15$
contains formulas of the form
$\sus q_1\,\kaz q_2\,\sus q_3\,\Phi(q_1,q_2,q_3)$,
where $\Phi$ is $\lip12$ and $q_i$ vary over $\bn$.

We allow codes $\rc\in\ccf_\om$
to substitute free variables over $\bn$.
If $\vpi:=\vpi(\rc_1,\dots,\rc_n)$ is an \dd\cL formula, 
and   $x\in\on$, 
\index{formula!$\vpi[x]$}%
\index{zzfix@$\vpi[x]$}%
then $\vpi[x]$ denotes the formula 
$\vpi(f_{\rc_1}(x),\dots,f_{\rc_n}(x))$, 
where all $f_{\rc_i}(x)$ are reals in $\bn$, of course.  

\bdf
[in $\rL$]
\lam{lfo}
We define a relation $T\fo\vpi$ between trees $T\in\pew$ 
and closed \dd\cL formulas in
$\bigcup_{k\ge2}(\lis1k\cup\lip1k)$.
Recall that $\tet=\omd$ (in $\rL$).
\ben
\Aenu
\itlb{lfo1}%
Let $\vpi(\rc_1,\dots,\rc_n)$ be a $\lis12$ formula
as in (1), and $\rc_1,\dots,\rc_n\in\ccf_\om$.
Let finally $T\in\pew$.
We define $T\fo\vpi$ iff
there exist codes                                
$\rc\in\ccf_\om$ and $\rd\in\ccf_{\tet}$
such that the following 
holds for all $x\in[T]$:
$$
\kaz m\,Q(f_\rd(x)\res m,f_{\rc}(x)\res m,
f_{\rc_1}(x)\res m,\dots,f_{\rc_n}(x)\res m).
$$    

\itlb{lfo2}%
If $\vpi$ is a closed $\lip1{k}$ formula, $k\ge 2$,
then $T\fo\vpi$ iff there is no tree $S\in\pew$   
such that $S\sq T$ and $S\fo\vpi^-$, where $\vpi^-$
is the result of canonical transformation of
$\neg\:\vpi$ to $\lis1k$ form.

\itlb{lfo3}%
If $\vpi:=\sus x\,\psi(x)$ is a closed  
$\lis1{k+1}$ formula, $k\ge 2$ 
($\psi$ being of type $\lip1k$),
then $T\fo\vpi$ iff there is a code $\rc\in\ccf_\om$   
such that $T\fo\psi(\rc)$.
\een
If $\vpi(p_1,\dots,p_n)$ is an 
\dd\cL formula then let
$$
\Fo\vpi=
\ens{
\ang{T,\rc_1,\dots,\rc_n}}{
T\in\pew
\land
\rc_i\in\ccf_\om
\land 
T\fo\vpi(\rc_1,\dots,\rc_n)
}.
$$
In particular if $\vpi$ is closed then 
$\Fo\vpi=\ens{T\in\pew}{T\fo\vpi}$. 
We also define 
$\des\vpi=\Fo\vpi\cup\Fo{\vpi^-}$ in this case.
\edf

\bte
[in $\rL$]
\lam{deff}
If\/ $k\ge2$ and\/ $\vpi$ is a formula in\/ 
$\lis1k$, resp., $\lip1k$, 
then the set\/ 
$\Fo\vpi$ belongs to\/  $\ahs{k-1}$, 
resp., $\ahp{k-1}$.
\ete
\bpf
The proof goes on by induction on $k$.
We begin with $\lis12$ formulas.
We argue in the assumptions and notation of (1) 
above. 
According to definition~\ref{lfo}\ref{lfo1}, 
the existence quantifiers over $\rc$ and $\rd$ are 
in line with the $\is{}1$ definability, but 
we have to prove that the set 
$$
\bay{l}      
W=\{\ang{\rd,\rc,\rc_1,\dots,\rc_n,T,m}\in 
\ccf_{\tet}\ti(\ccf_\om){}^{n+1}\ti\pew\ti\om:
\hspace*{8ex}\\[0.5ex]
\,\hspace*{17ex}
\kaz x\in[T]\,Q(f_\rd(x)\res m,f_{\rc}(x)\res m,
f_{\rc_1}(x)\res m,\dots,f_{\rc_n}(x)\res m)\}
{}
\eay
$$ 
belongs to $\ahs1$.               
Recall that $Q$ is $\ahd0$ by (1).
It can also be mentioned that 
$\ccf_{\tet}\cup\ccf_\om\cup\pew\sq\hcd$, 
so that $W\sq\hcd$ anyway.

The hostile elements in the definition of $W$, which do 
not allow it to be $\ahs1$ straightaway, are the 
quantifier $\kaz x\in[T]$ in the second line, and 
the quantifier $\kaz x\in\on$ in (2) of Definition~\ref{dccf}. 
(As we argue in $\rL$, the upper index $\rL$ as in \ref{dccf} 
is removed.) 
But, $\on\in\hcd$ (under $\rV=\rL$), hence, as we don't care 
here about the choice of parameters in $\hcd$,\snos
{If we do care then the result holds too but by means 
of more thoroughful arguments.} 
we can pick up $\on$ as the extra parameter. 
The quantifier $\kaz x\in\on$ in (2) of \ref{dccf} then 
immediately becomes bounded, while the quantifier 
$\kaz x\in[T]\,(\dots x\dots)$ in the definition of $W$ 
changes to 
$\kaz x\in\on\,(x\in[T]\imp \dots x\dots)$, 
hence becomes bounded as well, and overall we get even 
$W\in \ahd0$, as required.
     
The induction steps are easy applications of
\ref{lfo}\ref{lfo2},\ref{lfo3}.
\epf

Recall that a number $\nn\ge2$ is fixed by  
Definition~\ref{vdp}.

\ble
[in $\rL$]
\lam{dens}
Let\/ $\vpi$ be a closed formula in\/ $\lis1k\cup\lip1k$, 
$k\ge2$. 
Then the set\/ $\des\vpi$ is dense in\/ $\pew$. 
If\/ $k<\nn$, then\/ 
$\des\vpi\cap\dP$ is dense in\/ $\dP$.
\ele
\bpf
The first claim is a simple application of
Definition~\ref{lfo}\ref{lfo2}. 
The second claim follows from the first one by lemmas  
\ref{deff} and \ref{jden}\ref{prop3}.
\epf

\parf{Invariance} 
\las{inva}

It happens that the relation $\fo$ is invariant under 
some natural transformations of wide trees. 
Here we prove the invariance.
{\ubf We still argue in $\rL$.}

\vyk{
Let $S,T\in\pew$. 
There is a canonical \dd\su preserving bijection
$h_{ST}$ of $\br S$ onto $\br T$  
(and between each $\brn Sn$ and $\brn Tn$), defined
as follows.   
First, $h_{ST}(\roo S)=\roo T$.
Now suppose that $h_{ST}(s)=t$, $s\in \brn Sn$ and
$t\in \brn Tn$.
The sets $U=\ens{\xi<\omi}{s\we\xi\in S}$ and
$V=\ens{\et<\omi}{t\we\et\in T}$ have cardinality
$\card U=\card V=\ali$; let
$U=\ens{\xi_\ga}{\ga<\omi}$ and 
$V=\ens{\et_\ga}{\ga<\omi}$ be the enumerations in
the increasing order.
Now if $u\in\brn S{n+1}$ and $s\su u$ then
$u(\lh(s))=$ \text{some} $\xi_\ga$, and there is a
unique $v\in\brn T{n+1}$ with $t\su v$ and
$v(\lh(t))=\et_\ga$; we let $h_{ST}(u)=v$.
}

Let $S\in\pew$. 
To define a canonical homeomorphism 
$h_{S}:[S]\onto \on$, assume that $x\in[S]$.
Let $k<\om$.
Then $x\res{m_k}\in\brn Tk$ for some (unique) number $m_k$.
The set
$\Xi(x,k)=\ens{\xi<\omi}{(x\res{m_k})\we\xi\in S}$
has cardinality
$\card{(\Xi(x,k))}=\ali$; let
$\Xi(x,k)=\ens{\xi_\ga}{\ga<\omi}$ be the enumeration in
the increasing order.
In particular, $x(m_k)=\xi_\ga$ for some (unique)
$\ga=\ga(x,k)$.
Define $y=h_S(x)\in\on$ by $y(k)=\ga(x,k)$, $\kaz k$.
The map $h_S$ is a required homeomorphism.

It follows that if $T\in\pew$ is another tree then
$h_{ST}={h_T}\obr\circ h_S$ (the superposition)
is a homeomorphism of $[S]$ onto $[T]$.
Moreover, in this case, if $U\sq S$ is a subtree 
then the according subtree 
$h_{ST}\app U=\ens{h_{ST}(x)\res m}{x\in[U]\land m<\om}\sq T$ 
satisfies $U\in\pew$ iff $h_{ST}\app U\in \pew$,
and $[h_{ST}\app U]=\ens{h_{ST}(x)}{x\in [U]}$.   

\ble
[in $\rL$]
\lam{stuv}
If\/ $S,T\in\pew$ and\/ $U\in\pew$, $U\sq S$, then\/ 
$V=h_{ST}\app U\in\pew$, $V\sq T$, and\/ 
$h_{UV}=h_{ST}\res[U]$.\qed
\ele              

If $\la\in\Ord$ and  $f:[S]\to\la^\om$ then a function 
$h_{ST}\app f=f\circ h_{ST}\obr:[T]\to\la^\om$ 
is defined by $(h_{ST}\app f)(x)=f(h_{ST}(x))$, 
equivalently, $(h_{ST}\app f)(h_{ST}(x))=f(x)$.
If $\rc,\rc'\in\ccf_\la$ then we symbolically write 
$\rc'\res T=h_{ST}\app(\rc\res S)$, in case the associated 
functions $f_\rc$ and $f_{\rc'}$ satisfy: 
$f_{\rc'}\res[T]=f_{ST}\app(f_\rc\res[S])$.

\ble
[in $\rL$]
\lam{exc}
If\/ $S,T\in\pew$, $\la\in\Ord$, and\/ $\rc\in\ccf_\la$ 
then there is a code\/ $\rc'\in\ccf_\la$ satisfying\/ 
$\rc'\res T=f_{ST}\app(\rc\res S)$
\ele
\bpf
The map $f=f_\rc\res[S]:[S]\to\la^\om$ is continuous, 
hence so is the transformed map $f'=h_{ST}\app f:[T]\to\la^\om$. 
Let $g:\on\to\la^\om$ be any continuous extension of $f'$, 
and let $\rc'=\koe{g}$.
\epf

Finally if $\vpi:=\vpi(\rc_1,\dots,\rc_n)$ is a  
\dd\cL formula, and 
$\vpi':=\vpi(\rc'_1,\dots,\rc'_n)$, where 
$\rc'_1,\dots,\rc'_n$ 
is another set of codes $\rc'_i\in\ccf_\om$, 
then we symbolically write 
$\vpi'\res T=h_{ST}\app(\vpi\res S)$, in case 
$\rc_i'\res T=h_{ST}\app(\rc_i\res S)$ holds 
for each $i=1,\dots,n$.    

\ble
[in $\rL$]
\lam{inv}
Let\/ $S,T\in\pew$ and let\/ 
$\vpi,\vpi'$ be closed formulas in\/ $\lis1k\cup\lip1k$, $k\ge2$, 
and finally\/ $\vpi'\res T=h_{ST}\app(\vpi\res S)$.
Then\/   
$S\fo\vpi$ iff\/ $T\fo \vpi'$.
\ele
\bpf
We argue by induction. 
Let $\vpi,\vpi'$ be $\lis12$, so that 
$\vpi:=\vpi(\rc_1,\dots,\rc_n)$ and 
$\vpi':=\vpi(\rc'_1,\dots,\rc'_n)$, where 
$\rc_1,\dots,\rc'_n,\rc'_1,\dots,\rc'_n$ 
are codes in $\ccf_\om$, and 
$$
\vpi(p_1,\dots,p_n)
\,:=\, 
\sus\chi\in\tso\,\sus q\in\bn\,\kaz m\,
Q(\chi\res m,q\res m,p_1\res m,\dots,p_n\res m) 
$$ 
is a formula as in (1) of Section~\ref{fs},
and $\rc_i'\res T=h_{ST}\app(\rc_i\res S)$ holds 
for each $i$. 

Assume that $S\fo\vpi$. 
Then by definition (Definition~\ref{lfo}\ref{lfo1}) 
there are codes $\rc\in\ccf_\om$ and $\rd\in\ccf_{\tet}$ 
such that\pagebreak[0]  
$$
\kaz x\in[S]\,\kaz m\,Q(f_\rd(x)\res m,f_{\rc}(x)\res m,
f_{\rc_1}(x)\res m,\dots,f_{\rc_n}(x)\res m).
$$ 
Pick, by Lemma~\ref{exc}, codes $\rc'\in\ccf_\om$ and 
$\rd'\in\ccf_{\tet}$ with $\rc'\res T=h_{ST}\app(\rc\res S)$ 
and $\rd'\res T=h_{ST}\app(\rd\res S)$. 
Then we obtain\pagebreak[1]  
$$
\kaz y\in[T]\,\kaz m\,Q(f_{\rd'}(y)\res m,f_{\rc'}(y)\res m,
f_{\rc'_1}(y)\res m,\dots,f_{\rc'_n}(y)\res m)\,,
$$ 
and hence the codes $\rc'$ and $\rd'$ witness $T\fo \vpi'$.

{\it Step $\lis1k\to\lip1k$.} 
Let $\vpi$ be a closed formula in $\lip1k$, so that $\vpi$ 
is $\psi^-$, where $\psi$ is $\lis1k$, and accordingly 
$\vpi'$ is $(\psi')^-$, 
$\psi'\res T=h_{ST}\app(\psi\res S)$.
Assuming that $S\fo\vpi$, prove that $T\fo\vpi'$.
Suppose to the contrary that $T\fo\vpi'$ fails. 
Then, by Definition~\ref{lfo}\ref{lfo2}, there is a tree 
$V\in \pew$, $V\sq T$, $V\fo\psi'$.
We let $U=h_{TS}\app V$, so that $U\in\pew$, $U\sq S$, 
$V=h_{st}\app U$. 
And, by the way, $h_{UV}=h_{ST}\res[U]$ by Lemma~\ref{stuv}, 
thus still $\psi'\res V=h_{UV}\app(\psi\res[U])$. 
It follows that $U\fo\psi$, by the inductive hypothesis, 
which contradicts to $S\fo\vpi$. 

{\it Step $\lip1k\to\lis1{k+1}$.} 
Let $\vpi$ be a closed formula in $\lis1{k+1}$, 
so that $\vpi$ 
is $\sus q\,\psi(q)$, where $\psi(q)$ is $\lip1k$, and accordingly 
$\vpi'$ is $\sus q\,\psi'(q)$, 
$\psi'\res T=h_{ST}\app(\psi\res S)$.
Assuming that $S\fo\vpi$, prove that $T\fo\vpi'$.
By Definition~\ref{lfo}\ref{lfo3}, 
there is a code $\rc\in\ccf_\om$ 
satisfying $S\fo\psi(\rc)$.
By Lemma~\ref{exc}, there exists a code $\rc'\in\ccf_\om$ 
such that $\rc'\res T=h_{ST}\app(\rc\res S)$.
Then $\psi'(\rc')\res T=h_{ST}\app(\psi(\rc)\res S)$. 
It follows that $T\fo\psi'(\rc')$, by the inductive hypothesis, 
hence $T\fo\vpi'$.
\epf

\bcor
\lam{tt'}
Let\/ $S,T\in\pew$ and let\/ $\vpi$ be a closed formula in\/ 
$\lis1k\cup\lip1k$, $k\ge2$, 
with no codes in\/ $\ccf_\om$ as parameters.
Then\/   
$S\fo\vpi$ iff\/ $T\fo \vpi$.\qed
\ecor

\parf{Forcing and truth} 
\las{ftru}

Recall that $\nn\ge2$ is fixed by Definition~\ref{vdp}. 

Moreover we'll assume that $\nn\ge3$, because we 
now focus on the proof of claim \ref{Tun4} 
of Theorem~\ref{Tun}, vacuous in the case $\nn=2$.

The last part of the proof of Theorem~\ref{Tun} 
will be the next theorem which connects 
the forcing relation $\fo$ with the truth in 
\dd\dP generic extensions. 
This will be the key ingredient of the proof of 
Theorem~\ref{Tun}\ref{Tun4}: 
we use the invariant relation $\fo$ to 
surprisingly approximate the 
forcing $\dP$, definitely non-invariant under 
the transformations considered in Section~\ref{inva}.

\bte
\lam{f-t}
Assume that\/ $2\le k<\nn$, $\vpi\in\rL$ is a    
closed formula in\/ $\lip1{k}\cup\lis1{k+1}$, and a set\/ 
$G\sq\dP$ is generic over\/ $\rL$. 
Then the sentence\/ $\vpi[\xx G]$ is true in\/ $\rL[G]$  
if and only if\/ $\sus T\in G\,(T\fo\vpi)$.
\ete
\bpf
{\ubf We argue in $\rL[G]$.} 
{\ubf Base of induction:} 
$\vpi$ is a closed $\lis12$ formula,
$$
\vpi 
\,:=\, 
\vpi(\rc_1,\dots,\rc_n)
\,:=\, 
\sus\chi\in\tet^\om\,\sus q\in\bn\,\kaz m\,
Q(\chi\res m,q\res m,\rc_1\res m,\dots,\rc_n\res m)\,,
$$
as in \ref{lfo}\ref{lfo1} and (1) of Section~\ref{fs}. 
Assume that $T\in G$ and $T\fo\vpi$. 
Then by Definition~\ref{lfo}\ref{lfo1} there are codes 
$\rc\in\ccf_\om\cap\rL$ and $\rd\in\ccf_\tet\cap\rL$ such that 
$$
\kaz m\,\kaz x\in[T]\cap\rL\:Q(f_\rd(x)\res m,f_{\rc}(x)\res m,
f_{\rc_1}(x)\res m,\dots,f_{\rc_n}(x)\res m).
$$
(Recall Remark~\ref{rccf} on the absoluteness of being a code 
in any $\ccf_\la$.)
However all functions $f_\rd,f_{\rc},f_{\rc_i}$ are continuous. 
It follows that the last displayed formula can be strengthened 
to
$$
\kaz x\in[T]\,\kaz m\:Q(f_\rd(x)\res m,f_{\rc}(x)\res m,
f_{\rc_1}(x)\res m,\dots,f_{\rc_n}(x)\res m).
$$
Therefore, as $\xx G\in[T]$ (because $T\in G$), we obtain 
$$
\kaz m\:Q(f_\rd(\xx G)\res m,f_{\rc}(\xx G)\res m,
f_{\rc_1}(\xx G)\res m,\dots,f_{\rc_n}(\xx G)\res m). 
$$
Thus elements $\chi=f_\rd(\xx G)$ and 
$q=f_\rc(\xx G)$ witness $\vpi[\xx G]$ to be true. 

To establish the inverse, suppose that $\vpi[\xx G]$ 
is true in $\rL[G]$, that is,
$$
\kaz m\:Q(\chi\res m,q\res m,
f_{\rc_1}(\xx G)\res m,\dots,f_{\rc_n}(\xx G)\res m)  
$$
true for some $\chi\in\tet^\om$ and $q\in\bn$ in $\rL[G]$. 
By Lemma~\ref{repd} there are codes $\rd\in\ccf_\tet\cap\rL$ 
and $\rc\in\ccf_\om\cap\rL$ such that $\chi=f_\rd(\xx G)$ and 
$q=f_\rc(\xx G)$. 
Thus there is a tree $T\in G$ which \dd\dP forces the formula
$$
\kaz m\:Q(f_\rd(\xx \uG)\res m,f_{\rc}(\xx \uG)\res m,
f_{\rc_1}(\xx \uG)\res m,\dots,f_{\rc_n}(\xx \uG)\res m)  
\eqno(*)
$$
over $\rL$. 
We claim that the codes $\rc$ and $\rd$ witness $T\fo\vpi$ 
as in \ref{lfo}\ref{lfo1}. 
Indeed otherwise there are $x\in[T]$ and $m$ such 
that\pagebreak[0] 
$$
\neg\;Q(f_\rd(x)\res m,f_{\rc}(x)\res m,
f_{\rc_1}(x)\res m,\dots,f_{\rc_n}(x)\res m).
\eqno(\dag)
$$    
But, the maps $f_\rd,f_{\rc},f_{\rc_i}$ are continuous. 
It follows that there is a string $u=x\res j$ for some $j$ 
such that $(\dag)$ holds for {\ubf all} $x\in [S]$, 
where $S=\req Tu\in\dP$.
But then clearly $T$ cannot 
\dd\dP force $(*)$ as $S$ forces the opposite.

{\ubf Step $\lis1k\imp\lip1{k}$, $k<\nn$.} 
Let $\vpi$ be a $\lip1k$ formula. 
By Lemma~\ref{dens}, there is a tree $T\in G$  
such that either $T\fo\vpi$ or $T\fo\vpi^-$.
Assume that $T\fo\vpi$; 
we have to prove that $\vpi[\xx G]$ is true. 
Suppose otherwise. 
Then $\vpi^-[\xx G]$ is true. 
By the inductive hypothesis, there is a tree  
$S\in G$ such that $S\fo\vpi^-$.
But the trees $S,T$ belong to the same generic set 
$G$, hence they are compatible, which leads to a 
contradiction with the assumption $T\fo\vpi$, 
according to Definition~\ref{lfo}\ref{lfo2}.
Now assume that $T\fo\vpi^-$. 
Then $\vpi^-[\xx G]$ is true by  
the inductive hypothesis, hence $\vpi[\xx G]$ is false. 
On the other hand, there is no tree 
$S\in G$ such that $S\fo\vpi^-$, just as above.

{\ubf Step $\lip1k\imp\lis1{k+1}$, $k<\nn$.} 
Let $\vpi$ be $\sus x\,\psi(x)$ where $\psi$ 
is $\lip1k$. 
Assume that $T\in G$ and $T\fo\vpi$. 
Then by Definition~\ref{lfo}\ref{lfo3}  
there is a code $\rc\in\ccf\cap\rL$ such that 
$T\fo\psi(\rc)$.
By the inductive hypothesis, 
the formula $\psi(\rc)[\xx G]$, that is,
$\psi[\xx G](f_\rc(\xx G))$, is true in $\rL[G]$. 
But then $\vpi[\xx G]$ is  true as well. 

Conversely assume that $\vpi[\xx G]$ is true. 
Then there is a real $y\in\rL[G]\cap\bn$ such that
$\psi[\xx G](y)$ is true.
By Lemma~\ref{repd}\ref{repd1},  
$y=f_\rc(\xx G)$ for a code $\rc\in\ccf_\om\cap\rL$. 
But then $\psi(\rc)[\xx G]$ is true in $\rL[G]$. 
By the inductive hypothesis, there is a tree $T\in G$  
satisfying $T\fo\psi(\rc)$. 
Then $T\fo\vpi$ as well.
\epf

\parf{The final argument} 
\las{finar}

\bpf[\ubf Theorem~\ref{Tun}, the main theorem]
We assert that any \dd\dP generic extension 
$\rL[G]=\rL[\xx G]$ 
satisfies conditions (i), \ref{Tun2}, \ref{Tun3}, \ref{Tun4} 
of the theorem.
Regarding (i), \ref{Tun2}, \ref{Tun3} see a review in 
the very end of Section~\ref{bex1}. 
Let's concentrate on \ref{Tun4}. 
Let $\Phi(j)$ be a parameter-free $\is1\nn$ formula.\snos
{The case when $\Phi$ has real parameters in $\rL$ can also  
be handled with some extra care.} 
Thus\pagebreak[0] 
$$
\Phi(j)\,:=\,
\sus r_1\,\kaz r_2\,\dots\,
\kaze 
r_n\,\susf m\:
R_j(r_1\res m,r_2\res m,\dots,r_\nn\res m)\,, 
$$ 
where $r_i$ are variables over $\bn$, 
$R_j\sq(\nse){}^\nn$, 
$R_j\in \rL$, 
and the map $j\mto R_j$ is arithmetically 
definable in $\rL$.
Applying Theorem~\ref{abst} in $\rL$ with 
$\tet=\omdl=\omi^{\rL[G]}$ and $M=\rL[G]$, 
we get relations $Q_j=Q_\tet(R_j)$, and closed
$\lis1\nn$ formulas\pagebreak[0] 
$$
\bay{l}
\vpi_j
\::=\:
\sus r_1\,\kaz r_2\,\dots\,
\susf r_{n-2}\,
\kaz(\exists)\,\chi\in\tet^\om\,
\kaz(\exists)\,q\,
\sus(\forall)\,m\:\hspace*{15ex}\\[0.5ex]
\hspace*{25ex}
Q_j(\chi\res m,q\res m,r_1\res m,r_2\res m,\dots,r_{n-2}\res m)
\;,
\eay
$$          
satisfying $\Phi(j)\eqv\vpi_j$, $\kaz j$, both in $\rL$ and 
in any \dd\dP generic extension $\rL[G]$ of $\rL$.
It follows, by Theorem~\ref{f-t}, that the set 
$X=\ens{j}{\Phi(j)^{\rL[G]}}$ (defined in $\rL[G]$) 
satisfies $X=\ens{j}{\sus T\in G\,(T\fo\vpi_j)}$. 
Furthermore, as the formulas $\vpi_j$ do not contain 
codes in $\ccf_\om$, it follows, by Corollary~\ref{tt'}, 
that $X=\ens{j}{T\fo\vpi_j}$, where  $T$ is any 
particular tree in 
$\pel$, one and he same for all $j$.
We conclude that $X\in\rL$, as required. 
\epf

\parf{A problem} 
\las{vop1}

It is a challenge to figure out what kind of models 
the method of the 
proof of Theorem~\ref{Tun} gives for cardinals bigger than 
$\ali$. 
For instance, let $\pew_{\omd}$ be the set of all trees 
$T\sq{\omd}\lom$ whose all branching nodes are 
\dd\omd branching nodes.
This is a non-Laver version of the Namba forcing;
the Namba forcing per se requires that in addition 
every node above the stem is a branching node.
The forcing $\pew_{\omd}$ (or an equivalent forcing)
is considered
\eg\ in \cite{buk76}, \cite[Section 28]{jechmill}, and
\cite[18.4]{cumm}.

Clearly $\pew_{\omd}$ adds a cofinal infinite sequence,
say $\vec a=\sis{\al_n}{n<\om}$, 
in $\omd$. 
On the other hand, if CH holds in the ground
universe then, essentially by Namba,
$\pew_{\omd}$ does not add new
reals, hence, does not collapse $\omi$.
(See \cite[Section 28]{jechmill} for a simple proof.)
Thus $\vec a\in\hka\la$ in the extension $\rV[\vec a]$,
where $\la=\omd^{\rV[\vec a]}>\omd^\rV$.
(Where $\rV$ is the ground set universe, as usual.)
It is then an interesting problem to check whether
there are results for the definability of $\vec a$
in $\hka\la$ similar to the results in \cite{abr2}
and those of this paper.

\parf{Acknowledgements} 
\las{ask}

Vladimir Kanovei acknowledges partial support of
grant RFBR 17-01-00705, and is thankful to 
the Erwin Schrodinger International
Institute for Mathematics and Physics
for their support during the December 2016 visit.
Vassily Lyubetsky acknowledges
partial support of grant RSF 14-50-00150. 


\bibliographystyle{plain}
\addcontentsline{toc}{subsection}{\hspace*{5.5ex}References}
{\small

\bibliography{36}
}%

\vyk{
\def\indexname{\large Index\thanks{\\ \normalsize\ubf\ \ 
not a part of the manuscript, added for 
the convenience of the process of refereeing.}}%
}

\def\indexname{Index\,\footnotemark{}}

\indexprologue{
\footnotetext{\bf~The index is
not a part of the text, but is added for 
the convenience of the the process of refereeing 
of the manuscript.}%
\addcontentsline{toc}{subsection}{\hspace*{5.5ex}Index}%
\vspace*{-4ex}%
}

\small\printindex

\end{document}